\documentclass[12pt]{amsart}
\usepackage{mathrsfs}
\usepackage[colorlinks]{hyperref}
\usepackage{xcolor}
\usepackage{color}
\usepackage[a4paper,asymmetric]{geometry}
\usepackage{mathscinet}
\usepackage{fullpage}
\usepackage{latexsym}
\usepackage{amsthm}
\usepackage{amssymb}
\usepackage{amsfonts}
\usepackage{amsmath}
\usepackage{verbatim}
\newtheorem{theorem}{Theorem}[section]
\newtheorem{lemma}[theorem]{Lemma}
\newtheorem{proposition}[theorem]{Proposition}
\newtheorem{corollary}[theorem]{Corollary}
\newtheorem{hypothesis}[theorem]{Hypothesis}
\newtheorem{assumption}[theorem]{Assumption}
\theoremstyle{definition}
\newtheorem{definition}[theorem]{Definition}
\newtheorem{example}[theorem]{Example}
\theoremstyle{remark}
\newtheorem{remark}[theorem]{Remark}
\numberwithin{equation}{section}
\def\d{{\rm d}}
\def\e{{\rm e}}

\def\R{\mathbb R}

\def\la{\left(}
\def\ra{\right)}

\def\l2{L^2_{\theta,\delta}}

\def\ve{\varepsilon}

\begin{document}
\title[Linear parabolic equation with Dirichlet white noise boundary conditions]
{Linear parabolic equation with Dirichlet  white noise boundary conditions}
\author{Ben Goldys}
\address{School of Mathematics and Statistics, The University of Sydney, Sydney 2006, Australia}
\email{Beniamin.Goldys@sydney.edu.au}
\author{Szymon Peszat}
\address{Institute of Mathematics, Jagiellonian University, {\L}ojasiewicza 6, 30-348 Krak\'ow, Poland}
\email{napeszat@cyf-kr.edu.pl}
\thanks{This  work was partially supported by the ARC Discovery Grant DP120101886.}
\thanks{The work of Szymon Peszat    was supported by Polish National Science Center grant  2017/25/B/ST1/02584.}
\begin{abstract}
We study inhomogeneous Dirichlet boundary value problems associated to  a linear parabolic equation $\frac{du}{dt}=Au$ with strongly elliptic operator $A$ on bounded and unbounded domains with white noise  boundary data. Our main assumption is that the heat kernel of the corresponding homogeneous problem enjoys the Gaussian type estimates taking into account the distance to the boundary. Under mild assumptions about the domain, we show that $A$ generates a $C_0$-semigroup in weighted  $L^p$-spaces where the weight is a proper power of the distance from the boundary. We also prove some smoothing properties and exponential stability of the semigroup. Finally, we reformulate the  Cauchy-Dirichlet problem with white noise boundary data as an evolution equation in the weighted space and prove the existence of Markovian   solutions.  
\end{abstract}
\maketitle
\tableofcontents
\section{Introduction}\label{intro}
The aim of this paper is to study the following  linear stochastic boundary value problem 
\begin{equation}\label{E11}
\left\{\begin{array}{lll}
\dfrac{\partial X}{\partial t}(t,x)=\mathcal{A} X(t,x)\,,&x\in\mathcal{O},&t>0\,,\\
&\\
X(t,x)= \dfrac{\partial W}{\partial t}(t,x)\,,&x\in\partial\mathcal{O},&t>0\,,\\
&\\
X(0,x)=X_0(x)&x\in\mathcal{O}\,.&
\end{array}\right.
\end{equation}
where $\mathcal O\subset\R^d$ is an open, possibly unobunded, domain, $W$ is a Wiener process taking values in a space of distributions on $\partial \mathcal{O}$, and $\mathcal A$ is a  second order, strongly elliptic operator in $\mathcal{O}$. Let us note that solutions to \eqref{E11} are Markovian if and only if $W$ is a process with independent increments. 
\par
There exists vast literature on the non-homogeneous Dirichlet boundary value problem for deterministic linear parabolic equations, for a classical exposition see the fundamental monograph \cite{Lions-Magenes} or more recent \cite{mclean}. Extension of the classical results to rough boundaries and rough boundary conditions is still a subject of ongoing research, see for example \cite{Lindemulder-Veraar} and references therein. In this paper we study equation \eqref{E11} in a relatively regular domain,  see Section \ref{formulation} for details, but the boundary condition $\frac{\partial W}{\partial t}$ can be very irregular, including space-time white noise. Apart from purely mathematical motivations, such an extension is important in non-equlibrium statistical mechanics and optimal control theory, see for example \cite{Fabri-Goldys}, \cite{Maslowski}, \cite{duncan} and a recent book \cite{munteanu}. 
\par
Stochastic equations with boundary noise were usually studied in the case of Neumann boundary conditions that are more tractable, see \cite{Freidlin-Sowers}, \cite{Sowers}, \cite{munteanu_p} and also aforementioned papers \cite{Maslowski} and \cite{duncan}. Much less is known about stochastic equations with Dirichlet boundary noise. Equation \eqref{E11} was proposed in the seminal work \cite{DaPrato-Zabczyk}, where it was shown that it has no $L^2(\mathcal O,\d x)$-valued solutions. In \cite{Alos-Bonacorsi1} and \cite{Alos-Bonacorsi2} solutions to a nonlinear equation in $\mathcal O=(0,+\infty)$ for $\mathcal A=\frac{\d ^2}{\d x^2}$ are studied and proved to have trajectories in $L^2\la (0,+\infty);\,x^{1+\theta}\d x\ra$. In \cite{bgpr} a similar approach is used to consider a very general formulation of the stochastic boundary value problem in multidimensional domains for a large class of elliptic operators $\mathcal A$ and distribution-valued Gaussian noises, see also \cite{Brzezniak-Peszat} for the case of stochastic wave equation. 
\par
In \cite{Alos-Bonacorsi1, Alos-Bonacorsi2} the problem was not stated  as an evolution equation in $L^2\la (0,+\infty);\,x^{1+\theta}\d x\ra$. Such a formulation was introduced and exploited in \cite{Fabri-Goldys}. The main ingredient was a result by Krylov \cite{Krylov}, who proved that Laplacian generates a strongly continuous analytic semigroup in the space $L^2\la\R^d_+;\rho^{1+\theta}(x)\,\d x\ra$, where $\rho$ stands for the distance of a point $x\in\R^d_+$ to the boundary. It seems that the method used in \cite{Krylov} to prove this  generation result does not extend to more general domains and more general elliptic operators. 
One of our goals in this paper is to show that equation \eqref{E11} can be reformulated  as a  stochastic evolution equation 
\begin{equation}\label{E15}
\d X = AX\d t + B\d W, \quad X(0)= X_0,
\end{equation}
on a state space $E=L^p\la\mathcal O,\rho^\theta(x)\,\d x\ra$ with an appropriately chosen operator $B$. The operator $A$ is an abstract realisation of $\mathcal A$ as a generator of the $C_0$-semigroup in $E$. This will ensure  the  Markov property of the solution and since $E$ is a function space, it will open the way to study nonlinear perturbations  of \eqref{E11}. 
\section{Formulation of the problem}\label{formulation}
 The boundary noise $W$ is a Wiener process taking values in a space of distributions on $\partial \mathcal{O}$. More precisely, we will assume that $W$  can be represented  as  a formal  series 
\begin{equation}\label{E12}
W(t,x)= \sum_{k}e_k(x)W_k(t),
\end{equation}
where $W_k$ are independent real-valued Wiener processes defined on a filtered probability space  $(\Omega, \mathfrak{F}, (\mathfrak{F}_t),\mathbb{P})$ and  $(e_k)$ is  a finite or infinite sequence of functions on $\partial \mathcal{O}$. In order to simplify the presentation we assume that $\{e_k\}\subset  L^2(\partial \mathcal{O},\d \textrm{s})$,  where $\textrm{s}$ is the surface measure and that 
$$
\sum_{k} \left(\int_{\partial \mathcal{O}}e_k(y)\psi(y)\d \textrm{s}(y)\right)^2<+\infty, \qquad \forall\, \psi \in L^2(\partial \mathcal{O},\d \textrm{s}).
$$
However, our framework can be easily adapted to the case where $e_k$ are distributions on $\partial \mathcal{O}$, see however Remark  \ref{R88}.  

\begin{remark}\label{R11}
{\rm Let us recall,  see e.g. \cite{DaPrato-Zabczyk2},  that there exists a Hilbert space  $H_W$ called the \emph{Reproducing Kernel Hilbert Space} of $W$ such that  
$$
W(t,x)= \sum_{k}\tilde e_k(x)\tilde W_k(t),
$$
where $\tilde W_k$ are independent real-valued Wiener processes defined on a filtered probability space  $(\Omega, \mathfrak{F}, (\mathfrak{F}_t),\mathbb{P})$ and  $(\tilde e_k)$ is  an orthonormal basis of $H_W$.  It can be shown that  $\textrm{linspan}\left\{ e_k\right\}$ is a dense subspace of $H_W$.   In a particular case of the so-called space white  noise $H_W=L^2(\partial \mathcal{O},\d \textrm{s})$.  
} 
\end{remark}

In Section \ref{SFormal mild solution} we derive  the concept of a  \emph{formal mild solution} to \eqref{E11}. Briefly it is given by the formula 
\begin{equation}\label{E13}
X(t)= S(t)X(0)+ \int_0^t\left(\lambda-A\right)S(t-s)D_\lambda\d W(s),\qquad t\ge 0, 
\end{equation}
where $S$ is the semigroup generated by the realization $A$ of $\mathcal{A}$ with homogeneous boundary conditions,  and $D_\lambda$ is the  \emph{Dirichlet map}. Let us recall that given $\lambda\ge 0$ and a function $\gamma$ on $\partial \mathcal{O}$, $u=D_\lambda \gamma$ is,  the possibly weak, see Section \ref{SDirichlet},  unique solution to the Poisson equation 
\begin{equation}\label{E14}
\mathcal{A}u(x)=\lambda u(x),\quad x\in \mathcal{O}, \qquad u(x)=\gamma(x), \quad x\in \partial \mathcal{O}. 
\end{equation}

To the best of our knowledge, equation \eqref{E11} with the solution defined by \eqref{E13} has been introduced in \cite{DaPrato-Zabczyk}.

\begin{remark}\label{R12}
{\rm  The stochastic integral appearing in \eqref{E13} is not well defined  in a space $L^p(\mathcal{O})$ but it does exist in a certain space $E$ such that $ L^p(\mathcal{O})\hookrightarrow E$.  In fact, see \cite{DaPrato-Zabczyk}, Example \ref{Ex13}, Popositions \ref{P71}, \ref{P72}, \ref{P713}, \ref{P718},  the solution to the problem on a bounded interval, or half line or half-space lives in a Sobolev space of negative order or on weighted $L^p(\mathcal{O}, w(x)\d x)$-space. }
\end{remark}

It turns out, see \cite{Alos-Bonacorsi1, Alos-Bonacorsi2, bgpr}, that  under mild assumptions on $\mathcal{O}$,  $\mathcal{A}$, and $W$,   the solution $X$  to \eqref{E11} is a smooth ($C^\infty$ in time and space variables) random field on $(0,+\infty)\times \mathcal{O}$ and that there is a $\kappa >0$ such that for $t>0$, 
$$
\mathbb{E}\left\vert X(t,x)\right\vert ^p \le C(t) \left( \textrm{dist}\left(x,\partial {\mathcal{O}}\right)\right)^{-\kappa}. 
$$ 
Therefore, $X(t)$ takes values in  $L^p(\mathcal{O},w)$-space  with an appropriate weight function $w$. 

One of our main goals is to show that the problem \eqref{E11} can be written equivalently  as the stochastic partial differential equation 
\begin{equation}\label{E15a}
\d X = AX\d t + B\d W, \quad X(0)= X_0,
\end{equation}
on an appropriately chosen  state space $E$ with  $B=(\lambda-A)D_\lambda $.  This ensures  the  Markov property of the solution and if $E$ is a function-valued space, it will enable us to study nonlinear perturbations  of \eqref{E11}.  We have to face, however,   the problem with the interpretation of $(\lambda -A)D_\lambda $. In fact  $A$,  as the generator of the heat semigroup with homogeneous Dirichlet  boundary condition,   is defined on regular functions vanishing on the boundary $\partial \mathcal{O}$,  whereas the restriction of $D_\lambda\frac{\partial W}{\partial t}$ to $\partial \mathcal{O}$ equals $\frac{\partial W}{\partial t}$! Therefore one needs to consider $A$ as the generator of the extension of the  semigroup $S$ on a suitable Sobolev space of negative-order (subspace of the space of distributions). The following example is taken from \cite{bgpr}. 
\begin{example}\label{Ex13}
Assume that $\mathcal{O}=(0,1)$,  $\mathcal{A}= \frac{\d ^2}{\d x^2}$, and $\lambda=0$. Then any function $\gamma \colon \partial \mathcal{O}=\{0,1\}\mapsto \mathbb{R}$ can be identified  with a pair $(\gamma_0,\gamma_1)\in \mathbb{R}^2$. We have 
$$
D_0(\gamma_0,\gamma_1)(x)=\gamma_0+ (\gamma_1-\gamma_0)x, \qquad x\in (0,1),
$$
and 
$$
AD_0(\gamma_0,\gamma_1)= \gamma_0\delta_0'-\gamma_1\delta_1',
$$
where $\delta'_a$ is the derivative of the Dirac delta distribution et $a$, and $A$ is the generator of the heat semigroup considered, for example, on the Sobolev space $W^{2,-2}(0,1)$. 
\end{example}

\begin{hypothesis}\label{H14}
There are $\lambda \ge 0$, $p>1$ and $s_0\ge 0$ such that the Dirichlet map $D_\lambda$ is a well defined bounded linear operator acting from $\textrm{linspan}\left\{ e_k\right\}$ into the Sobolev space $W^{-s_0,p}(\mathcal{O})$. 
\end{hypothesis}

\begin{hypothesis}\label{H15}
Operator $\mathcal{A}$ with homogeneous Dirichlet boundary conditions generates an analytic $C_0$-semigroup $S$ on each $W^{s,p}(\mathcal{O})$-spaces. For all $s,s'\in \mathbb{R}$, $p>1$ and $t>0$, $S(t)\colon W^{s,p}(\mathcal{O})\mapsto W^{s',p}(\mathcal{O})$. Moreover, if $A_{s,p}$\footnote{Later we will skip the subscripts $s$ and $p$ and we will write $A$ instead of $A_{s,p}$.} denotes the generator of $S$ on $W^{s,p}(\mathcal{O})$, then we assume that there is an $s_1$ such that $W^{-s_0,p}(\mathcal{O})  \hookrightarrow  D(A_{-s_1,p})$.
\end{hypothesis}
\begin{remark}\label{R26}
{\rm It is well known that Hypotheses \ref{H14} and \ref{H15} hold in a number of cases. By Theorem 4.10 in \cite{mclean} if $\mathcal  O$ is a bounded Lipschitz domain and the operator $A$ has Lipschitz coefficients, then $D_0\colon H^{1/2}(\partial \mathcal{O})\to H^1(\mathcal{O})$ is well defined and bounded.  In that case it is enough to assume that $\textrm{linspan}\left\{ e_k\right\}\subset H^{1/2}(\partial \mathcal{O})$.

If $\mathcal  O$ is a bounded $C^\infty$ domain and the operator $A$ has $ C^\infty$ coefficients, then  
$$
D_0\colon H^{-s-\frac{3}{2}}(\partial\mathcal {O})\to H^{-s}(\mathcal{O})
$$
is well defined and bounded for any $s\ge 0$, see Sections 6 and 7 in Chapter 2 of \cite{Lions-Magenes}. In particular, if $s\le -\,\frac{3}{2}$ then  $\textrm{linspan}\left\{ e_k\right\}\subset H^{-s-\frac{3}{2}}(\partial \mathcal {O})\subset L^2(\mathcal{O})$. \\

Very general conditions given in terms of capacities of $\mathcal{O}$ can be found in Chapter 15.7 of \cite{mazya}. }
\end{remark}

Hypotheses \ref{H14} and \ref{H15} enable us to reformulate problem \eqref{E13} into problem \eqref{E15a}  considered on the state space  $W^{-s_1,p}(\mathcal{O})$. In fact the map 
$$
B= \left(\lambda -A\right)D_\lambda :=  \left(\lambda -A_{-s_1,p}\right)D_\lambda
$$
is a bounded linear operator from $\textrm{linspan}\left\{ e_k\right\}$  into $W^{-s_1,p}(\mathcal{O})$ and $A= A_{-s_1,p}$ generates a $C_0$-semigroup $S= S_{-s_1,p}$ on $W^{-s_1,p}(\mathcal{O})$. Therefore, by our Proposition \ref{P61}  we have the following result.
\begin{theorem}\label{T16}
Under Hypotheses \ref{H14} and \ref{H15}, problem  \eqref{E15a} has the mild solution  solution 
\begin{equation}\label{E16}
X(t)= S(t)X_0+\int_0^t S(t-s)B\d W(s)
\end{equation}
in $W^{- s_1,p}(\mathcal{O})$-space if and only if 
\begin{equation}\label{E17}
\int_{\mathcal{O}} \left[ \sum_k  \int_0^T\left(\left( I-\Delta\right)^{-s_1/2} S(t) Be_k\right)^2(x)\d t \right]^{p/2}\d x <+\infty 
\end{equation}
for a certain or equivalently for any $T\in (0,+\infty)$. Moreover, if there is an $\alpha >0$ such that 
$$
\int_{\mathcal{O}} \left[ \sum_k  \int_0^Tt^{-\alpha}\left(\left( I-\Delta\right)^{-s_1/2} S(t) Be_k\right)^2(x)\d t \right]^{p/2}\d x <+\infty 
$$
then the mild solution has continuous trajectories\footnote{In fact H\"older continuous with arbitrary exponent $<\alpha/2$.}  in $W^{- s_1,p}(\mathcal{O})$, 
\end{theorem}

\begin{remark}
{\rm In Section \ref{SRadonifying} we will show that  condition  \eqref{E17}  guarantees that for any $t\ge 0$, stochastic integral $\int_0^t S(t-s)B\d W(s)$ is well defined in $W^{- s_1,p}(\mathcal{O})$. Note that, if $p=2$,  than $W^{- s_1,p}(\mathcal{O})$ is Hilbert space, and \eqref{E17} can be equivalently written as 
$$
\int_0^ T \|S(t)B\|^2_{L_{(HS)}(H_W, W^{- s_1,2}(\mathcal{O}))}\d t<+\infty, 
$$
where$\|\cdot\|_{L_{(HS)}(H_W, W^{- s_1,2}(\mathcal{O}))}$ is the Hilbert--Schmidt norm and $H_W$ is the Reproducing Kernel Hilbert Space of $W$, see Remark \ref{R11}. 

Since 
$$
\left(\lambda-A_{0,p}\right) S_{0,p}(s)D_\lambda =  S_{-s_1,p}(s)\left(\lambda-A_{-s_1,p}\right)D_\lambda
$$
the formal  mild solution and the mild solution defined by \eqref{E16}  coincide. }
\end{remark}

In order to obtain the function-valued solutions we need the following assumption. 
\begin{hypothesis}\label{H18}
The semigroup $S$ can be extended  to a $C_0$-semigroup on the weighted space  $L^p_{\theta,\delta}:= L^p(\mathcal{O},w_{\theta,\delta}(x)\d x)$, where 
\begin{equation}\label{E18}
w_{\theta,\delta}(x)=\min\left\{\text{\rm dist}\left(x, \partial \mathcal{O}\right)^\theta,  \left(1+\vert x \vert  ^2\right)^{-\delta}\right\}, 
\end{equation}
$p\in (1,+\infty)$, $\theta \in [0,2p-1)$ and $\delta \ge 0$. 
\end{hypothesis}

In Section \ref{SSemigroup} we will show that Hypothesis \ref{H18} is fulfilled under very mild assumptions on $\mathcal{O}$ and $\mathcal{A}$.  

Under the above three hypotheses the operator $B$ acts from $\textrm{linspan}\left\{ e_k\right\}$ into $W^{-s_0,p}(\mathcal{O})$. For any $t>0$, the $C_0$-semigroup $S(t)$ on $L^p_{\theta,\delta}$ has a unique continuous extension 
$$
S(t)\colon W^{-s_0,p}(\mathcal{O})\mapsto W^{0,p}(\mathcal{O})= L^p(\mathcal{O}, \d x)\hookrightarrow L^p_{\theta,\delta} .
$$
Therefore, as a consequence of our Proposition \ref{P61} and the classical theory of SPDEs (see e.g. \cite{DaPrato-Zabczyk})  we have the following general result. 
\begin{theorem} \label{T19}
Assume Hypotheses \ref{H14}, \ref{H15}, \ref{H18}. Problem  \eqref{E15a} has the mild solution  solution in $L^p_{\theta,\delta}$ if and only if 
\begin{equation}\label{E19}
\mathcal{J}_T(\{e_k\},p,\theta, \delta):= \int_{\mathcal{O}}\left[  \sum_{k} \int_0^T \left( S(t)Be_k\right)^2(x) \d t \right] ^{p/2} w_{\theta,\delta}(x)\d x <+\infty
\end{equation}
for a certain or equivalently for any $T\in (0,+\infty)$. Moreover, \eqref{E19} guarantees that problem \eqref{E15a} equivalently \eqref{E11}, defines a Markov family on the  state  space $L^p_{\theta,\delta}$.  If for a certain $\alpha >0$, 
\begin{equation}\label{E110}
\mathcal{J}_{T,\alpha} (\{e_k\},p,\theta, \delta):= \int_{\mathcal{O}}\left[  \sum_{k} \int_0^T t^{-\alpha}\left( S(t)Be_k\right)^2(x) \d t \right] ^{p/2} w_{\theta,\delta}(x)\d x <+\infty,
\end{equation}
then the mild solution has continuous trajectories in $L^p_{\theta,\delta}$.

 Finally,  the existence of an invariant measure is equivalent to the integrability condition 
\begin{equation}\label{E111}
\mathcal{J}_{+\infty}(\{e_k\},p,\theta, \delta):= \int_{\mathcal{O}} \left[ \sum_k \int_0^{+\infty}  \left( S(t)Be_k\right) ^2\d t \right] ^{p/2} w_{\theta,\delta}(x)\d x<+\infty. 
\end{equation}
\end{theorem}

\begin{remark}
{\rm If the semigroup $S$ is exponentially stable, i.e.  for a certain $\alpha >0$,
$$
\|S(t)\|_{L(L^p_{\theta,\delta}, L^p_{\theta,\delta})}\le C\e^{-\alpha t}, \qquad t\ge 0,
$$
 then condition  \eqref{E111} follows from \eqref{E19}.  In  Theorem \ref{T42}, we will show that the semigroup $S$ is exponential stable on $L^p_{\theta,\delta}$ if it is exponentially stable on $L^p_{0,\delta}$. Obviously if the domain $\mathcal{O}$ is bounded then for all $p$, $\theta$ and $\delta$, the spaces $L^p_{\theta,\delta}$ and $L^p_{\theta,0}$ are equivalent. Therefore, if $\mathcal{O}$ is bounded then we can always take $\delta=0$. Note that if $\mathcal{O}$ is bounded and $\mathcal{A}$ equals Laplace operator $\Delta$, then the corresponding semigroup is exponentially stable on $L^p_{0,0}$ and consequently on $L^p_{\theta,0}$ for any $p>1$. }
\end{remark}
 \begin{remark}
 {\rm Assume \eqref{E19}, Then for any $X_0\in L^p_{\theta,\delta}$ and for any $t>0$, $X(t)$ is a gaussian element in $L^p_{\theta,\delta}$. Therefore, by the Fernique theorem there is a $\beta >0$ such that 
 $$
 \mathbb{E}\, \exp\left\{ \beta \left\vert  X(t)\right\vert ^2_{L^p_{\theta,\delta}}\right\}<+\infty. 
 $$
 If \eqref{E110} is satisfied for an $\alpha >0$, then for any $T\in (0,+\infty)$, and for any $X_0\in L^p_{\theta,\delta}$, $X(\cdot )$ is a gaussian random element in $C([0,T]; L^p_{\theta,\delta})$. Thus there is a $\beta >0$ such that 
 $$
 \mathbb{E}\, \exp\left\{ \beta \sup_{t\in [0,T]}\left\vert  X(t)\right\vert ^2_{L^p_{\theta,\delta}}\right\}<+\infty. 
 $$
 }
 \end{remark}
Our framework enables us to study nonlinear problems. 
\begin{theorem} 
Assume \eqref{E19}, and Hypotheses \ref{H14}, \ref{H15}, \ref{H18}. Then for any  Lipschitz  continuous function $f\colon \mathbb{R}\mapsto \mathbb{R}$, and any $X_0\in L^p_{\theta,\delta}$, the boundary problem 
$$
\left\{\begin{array}{lll}
\dfrac{\partial X}{\partial t}(t,x)=\mathcal{A} X(t,x) + f(X(t,x))\,,&x\in\mathcal{O},&t>0\,,\\
&\\
X(t,x)= \dfrac{\partial W}{\partial t}(t,x)\,,&x\in\partial\mathcal{O},&t>0\,,\\
&\\
X(0,x)=X_0(x),&x\in\mathcal{O},\,&
\end{array}\right.
$$
has a unique solution in $L^p_{\theta,\delta}$, and 
$$
X(t)= S(t)X_0+ \int_0^t S(t-s)F(X(s))\d s + \int_0^t S(t-s)B\d W(s), 
$$
where $F(X(t))(x)= f(X(t,x))$. Finally, if the semigroup is exponentially stable with exponent $L$  and the Lipschitz constant of $f$ strictly less than $L$, then there is a unique invariant measure on $L^p_{\theta,\delta}$ for the nonlinear problem. 
\end{theorem}

The paper is organized as follows. In the next section we will heuristically derive the concept of the formal mild solution. Then, in Section \ref{SSemigroup}, we will show that Hypothesis \ref{H18} about  the $C_0$-property  of the semigroup on weighted $L^p$-spaces, is fulfilled  under rather mild  assumptions. The main difficulty is that the weight of the above form is not an $\mathcal {A}^*$-excessive function. Therefore the semigroup is not of contraction type. In Section \ref{SStability} we will study properties of the semigroup on weighted spaces. In our opinion the results of Sections \ref{SSemigroup} and \ref{SStability} are of independent interest.  

In Section \ref{SDirichlet} we will derive some useful point estimates for $S(t)Be$ for $e\in L^2(\partial \mathcal{O},d\mathrm{s})$.  In Section \ref{SRadonifying} we outline  the concept of stochastic integration in $L^p$-spaces.   Section \ref{SExamples} is devoted  to particular examples. We relay on estimates established in Section \ref{SDirichlet} and on results from the previous section. 

\section{Formal mild solution}\label{SFormal mild solution}

In our derivation of the concept of the \emph{formal mild solution to} \eqref{E11} we follow \cite{DaPrato-Zabczyk}. Recall that $D_\lambda$ denotes the Dirichlet map (see  \eqref{E14}).  Assume Hypothesis \ref{H14} and \ref{H15}.

Assume temporally that the boundary perturbation is of the form 
$$
W(t,x)=\sum_k e_k(x)\beta _k(t), 
$$
where series is finite, $e_k$ are functions or distributions on $\partial \mathcal{O}$, and $\beta_k\in C^1([0,+\infty))$.   We assume that for any $k$, $e_k$ belongs to the domain of the Dirichlet map $D_\lambda$ and that $D_\lambda e_k \in W^{-s_0, p}(\mathcal{O})$.

Note that if   $X$ is a solution to \eqref{E11} with $W$ as above, then 
$$
Y(t,x):= X(t,x)-D_\lambda\dfrac{\partial W}{\partial t}(t,x)
$$
satisfies  the homogeneous Dirichlet boundary conditions. Moreover, at least formally,  for $t>0$ and $x\in \mathcal{O}$ we have 
\begin{align*}
\frac{\partial Y}{\partial t}(t,x)&= AX(t,x)-\dfrac{\partial }{\partial t}D_\lambda \dfrac{\partial W}{\partial t}(t,x)\\
&= AY(t,x)+ \lambda D_\lambda \dfrac{\partial W}{\partial t}(t,x)-\dfrac{\partial }{\partial t}D_\lambda\dfrac{\partial W}{\partial t}(t,x). 
\end{align*}
Therefore 
\begin{align*}
Y(t,x)&= S(t)Y(0,x) +\int_0^t S(t-s)\left[ \lambda D_\lambda \dfrac{\partial W}{\partial s}(s,x)-\dfrac{\partial }{\partial s}D_\lambda \dfrac{\partial W}{\partial s}(s,x)\right]\d s\\
&= S(t)Y(0,x)+\int_0^t S(t-s)\lambda D_\lambda \dfrac{\partial W}{\partial s}(s,x)\d s \\
&\qquad - \left[ D_\lambda \dfrac{\partial W}{\partial t}(t,x)- S(t)D_\lambda \dfrac{\partial W}{\partial t}(0,x) + \int_0^t AS(t-s)D_\lambda \dfrac{\partial W}{\partial s}(s,x)\d s \right]\\
&= S(t)X(0,x) -D_\lambda \dfrac{\partial W}{\partial t}(t,x) +\int_0^t \left(\lambda -A\right)S(t-s)D_\lambda \dfrac{\partial W}{\partial s}(s,x)\d s. 
\end{align*}
Hence we infer that 
\begin{equation}\label{E21}
X(t)= S(t)X(0)+ \int_0^t \left(\lambda-A\right)S(t-s)D_\lambda  \d W(s).  
\end{equation}

Let us recall that the space $W^{-s_1,p}(\mathcal{O})$ appears in Hypotheses \ref{H15}.  Note that for any $k$ the stochastic integral 
$$
\int_0^t \left(\lambda -A\right) S(t-s)D_\lambda e_k \d W_k(s),\qquad t\ge 0, 
$$
takes valued in $W^{-s_1,p}(\mathcal{O})$. 
\begin{definition}
Let $X(0)\in W^{-s_1,p}(\mathcal{O})$, and let  the noise $W$ in $\eqref{E11}$  have  the form   \eqref{E12}. If  the series 
$$
\sum_{k} \int_0^t \left(\lambda -A\right) S(t-s)D_\lambda e_k \d W_k(s)=: \int_0^t \left(\lambda -A\right)S(t-s)D_\lambda \d  W(s)
$$
converges in $L^p(\Omega,\mathfrak{F},\mathbb{P};W^{-s_1,p}(\mathcal{O}))$, then  we call the proces defined by  formula \eqref{E21} the \emph{formal mild solution} to \eqref{E11} in $W^{-s_1,p}(\mathcal{O})$. 
\end{definition}

\section{Semigroup in weighted spaces}\label{SSemigroup}
Let $\mathcal{O}\subset\mathbb{R}^d$, $\mathcal{O}\not = \mathbb{R}^d$,  be an open connected domain.  From now on the following two assumptions will be satisfied.  

\begin{assumption}\label{A31}
We will assume that $\mathcal{O}$ is a $C^{1,\alpha}$-domain with $\alpha \in (0,1)$,   satisfying the connected line condition, see e.g. \cite{Chi-Kim-Park} for a  precise definition. Let us recall here that the connected line condition holds in many important cases including: 
\begin{itemize}
\item 
bounded $C^{1,\alpha}$ domain,
\item
graph above $C^{1,\alpha}$ function,
\item
$\mathcal{O}=\mathbb{R}^d_+$ or 
$$
\mathcal{O}=\left\{\left( x_i\right)\in\mathbb{R}^{d+1}\colon  a<x_{d+1}<b\right\}. 
$$
\end{itemize}
\end{assumption}

Let us consider  the following  second order differential operator
$$
\mathcal{A}\phi(x)=\sum_{i,j=1}^d\frac{\partial}{\partial x_i}\left( a_{ij}(x)\frac{\partial\phi}{\partial x_j}(x)\right)+\sum_{i=1}^d \mu^i(x)\frac{\partial\phi}{\partial x_i}(x).
$$
\begin{assumption}\label{A32}
We assume that  the homogeneous Dirichlet boundary problem  
\begin{equation}\label{E31}
\left\{\begin{array}{lll}
\dfrac{\partial u}{\partial t}(t,x)=\mathcal A u(t,x)\,,&x\in\mathcal{O},&t>0\,,\\
&\\
u(t,x)=0\,,&x\in\partial\mathcal{O},&t>0\,,\\
&\\
u(0,x)=f(x)&x\in\mathcal{O}\,,&
\end{array}\right.
\end{equation}
generates  a $C_0$-semigroup $(S(t))$ on $L^2(\mathcal{O}, \d x)$. The  generator of this semigroup will be denoted by $A$. Next, we assume that the semigroup can be represented by a Green kernel $G$, 
\begin{equation}\label{E32}
S(t)\psi(x)=\int_{\mathcal{O}}G(t,x,y)\psi(y)\d y,\quad x\in\mathcal{O}.
\end{equation}
Finally we assume that there exists a constant $\lambda>0$ such that 
$$
\lambda \vert h\vert ^2\le\langle a(x)a^T(x)h,h\rangle\le\lambda^{-1}\vert h\vert ^2,\quad x,h\in\mathbb{R}^d, 
$$
\begin{equation}\label{E33}
G(t,x,y)\le Cm_t(y)g_{ct}(x-y), \qquad t\le 1, \qquad x,y\in \mathcal{O},
\end{equation}
and
\begin{equation}\label{E34}
\left|\nabla_xG(t,x,y)\right|\le C\frac{m_t(y)}{\sqrt{t}}g_{ct}(x-y), \qquad t\le 1, \qquad x,y\in \mathcal{O},
\end{equation}
where 
$$
m_t(z):=\min\left\{ 1, \frac{\rho(z)}{\sqrt{t}}\right\}, \qquad \rho(z):=\text{\rm dist}\left(z, \partial \mathcal{O}\right) 
$$
and 
$$
g_t(z)= \left(2\pi t\right)^{-\frac{d}{2}}\e^{-\frac{\vert z\vert ^2}{2t}}. 
$$
\end{assumption}
\begin{remark}\label{R33}
{\rm Assumption \ref{A32}  is fulfilled if Assumption  \ref{A31} holds, the operator $\mathcal{A}$ is uniformly elliptic, the coefficients $a_{ij}$ are Dini continuous,  and $\mu^i$ are sign measures of the parabolic Kato class. In general $a_{i,j}$ and $\mu^i$  may  depend on $t$ and $x$ variables. For more details see \cite{Chi-Kim-Park}. In fact in \cite{Chi-Kim-Park} the following stronger estimate has beed obtained 
$$
G(t,x,y)\le Cm_t(x)m_t(y)g_{ct}(x-y), \qquad t\le 1, \qquad x,y\in \mathcal{O}. 
$$
}
\end{remark}

In the main theorem of this section we require the following assumption
\begin{assumption}\label{Axx}
For any $c>0$ and $\alpha \in (-1,0)$ there is a constant $C<+\infty$ such that 
$$
\sup_{x\in \mathcal{O}}\int_{\mathcal{O}}\rho^\alpha(y)g_{ct}(x-y)\d y \le Ct^{\frac{\alpha}{2}}, \qquad \forall\, t\in (0,1].  
$$
\end{assumption}
A proof of the following lemma is postponed to Appendix \ref{AppL}. 
\begin{lemma}\label{L45}
Assumption \ref{Axx} is satisfied if $\mathcal{O}$ is a half space or if $\mathcal{O}$  is a bounded $C^{1,\alpha}$-domain. 
\end{lemma}

\bigskip 
Recall  that the family of weights $w_{\theta,\delta}$, $\theta\ge 0$, $\delta \ge 0$,  were introduced in \eqref{E18}.  We will use the notations
$$
L^p_{\theta,\delta}:=L^p\left(\mathcal{O}, \mathcal{B}(\mathcal{O}),w_{\theta,\delta}(x)\d x\right), \qquad  L^p:=L^p_{0,0}= L^p(\mathcal{O},\d x),  \qquad p\ge 1, \  \theta, \delta\ge 0. 
$$

Let  $S=(S(t))$ be  the  $C_0$ semigroup on $L^2$ corresponding to \eqref{E31}. By Assumption  \ref{A32},  for each $t>0$, $S(t)$ is  defined by \eqref{E33} at least on  compactly supported functions $\psi$. 

The main result of this section is the following theorem.  Its proof is  given in  Section \ref{SProof}. 
\begin{theorem} \label{T34}
Let $p\in [1,+\infty)$, $\theta \in [0,2p-1)$ and $\delta\ge 0$.  Under Assumptions \ref{A31},  \ref{A32}, and \ref{Axx} we have: 
\begin{itemize}
\item[(i)]
For each $t$, $S(t)$ defined on compactly supported functions by \eqref{E32} has a unique extension to a bounded linear operator, denoted still by $S(t)$, acting from $L^p_{\theta,\delta}$ into $L^p_{\theta,\delta}$. Moreover, $S=(S(t))$  forms a   $C_0$-semigroup on $L^p_{\theta,\delta}$. 
\item[(ii)]  There exists a constant $C>0$ such that for all $t\in(0,1]$ and $\psi\in L^p_{\theta,\delta}$,  $S(t)\psi(x)$ is differentiable for each $x\in \mathcal{O}$ and 
$$
\left\vert  \frac{\partial}{\partial x_i} S(t)\psi  \right\vert _{L^p_{\theta,\delta}}\le \frac{C}{\sqrt{t}}|\psi|_{L^p_{\theta,\delta}},\qquad i=1,\ldots, d. 
$$
\end{itemize}
\end{theorem}
\begin{remark}
{\rm If $\theta >2p-1$, then  $L^p_\theta$ contains   functions $f$ with growth  $\rho^{-2}(y)$ at vicinity of some point of $\partial \mathcal{O}$.  On the other hand, the integral  $\int_{\mathcal{O}}G(t,x,y)f(y)\d y$  does not converge as $G(t,x,y)$ decays only at rate $\rho(y)$ at the boundary. Therefore, for $t>0$, $S(t)$ cannot be extended to $L^p_\theta$. } 
\end{remark}
\begin{remark}
{\rm For $0\le \theta <p$ we are able to  show the $C_0$-property and gradient estimates without  Assumption \ref{Axx}, for details see Appendix \ref{AppT}.}
\end{remark} 
\subsection{Preliminaries} 
Let $w^i\colon \mathcal{O}\mapsto (0,+\infty)$, $i=1,2$,  be measurable weights. Let 
$$
\mathcal{L}_i^p := L^p\left(\mathcal{O}, \mathcal{B}(\mathcal{O}), w^i(x)\d x\right),\qquad i=1,2\,,
$$
and let 
$$
\mathcal{L}^p := L^p(\mathcal{O}, \mathcal{B}(\mathcal{O}), w(x)\,\d x),
$$
where $w (x)=\min\{w^1(x),w^2(x)\}$. We will need the following elementary result. 
\begin{lemma}\label{L34}
Assume that $T$ is a bounded linear operator from $\mathcal{L}_i^p$ to $\mathcal{L}_i^p$ for $i=1,2$. Then it is bounded from $\mathcal{L}^p$ to $\mathcal{L}^p$ and the operator norm satisfies the estimate 
$$
\|T\|_{L(\mathcal{L}^p)}\le N:=2^{(p-1)/p}\max\left\{ \|T\|_{L(\mathcal{L}_1^p)},\|T\|_{L(\mathcal{L}_2^p)} \right\}.
$$
\end{lemma}
\begin{proof} Let $\psi\in \mathcal{L}_1^p\cap \mathcal{L}_2^p$ and 
$$
\mathcal{D}:= \{x\in \mathcal{O}\colon w_1(x)<w_2(x)\}, \qquad \mathcal{D}^c:= \mathcal{O}\setminus \mathcal{D}. 
$$
Then 
\begin{align*}
\int_{\mathcal O}\left\vert T\psi(x)\right\vert ^pw(x) \d x &= \int_{\mathcal O}\left\vert T(\chi_{\mathcal{D}}\psi)(x)
+T(\chi_{\mathcal{D}^c}\psi)(x)\right\vert ^pw(x) \d x\\
&\le 2^{p-1} \left[ \int_{\mathcal O}\left\vert T(\chi_{\mathcal{D}}\psi)(x)\right\vert^p w(x) \d x+ \int_{\mathcal O}\left\vert
T(\chi_{\mathcal {D}^c}\psi)(x)\right\vert ^pw(x) \d x\right]\\
&\le 2^{p-1} \left[ \int_{\mathcal O}\left\vert T(\chi_{\mathcal{D}}\psi)(x)\right\vert^p w^1(x) \d x+ \int_{\mathcal  O}\left\vert
T(\chi_{\mathcal {D}^c}\psi)(x)\right\vert ^pw^2(x) \d x\right]\\
&\le N^p  \left[ \int_{\mathcal O}\left\vert \chi_{\mathcal{D}}(x)\psi(x)\right\vert^p w^1(x) \d x+ \int_{\mathcal O}\left\vert
\chi_{\mathcal{D}^c}(x)\psi(x)\right\vert ^pw^2(x) \d x\right]. 
\end{align*}
Since 
$$
  \int_{\mathcal O}\left\vert \chi_{\mathcal {D}}(x)\psi(x)\right\vert^p w^1(x) \d x+ \int_{\mathcal  O}\left\vert
\chi_{\mathcal {D}^c}(x)\psi(x)\right\vert ^pw^2(x) \d x= \vert \psi\vert ^p_{\mathcal{L}^p},
$$
we have the desired conclusion. 
\end{proof}

\subsection{Proof of Theorem \ref{T34}}\label{SProof} 
Let 
$$
L^p_\theta:=L^p\left(\mathcal{O},\rho^\theta(x)\d x\right)\quad \text{and}\quad \mathcal{L}^p_\delta:= L^p\left(\mathcal{O}, \left(1+\vert x\vert ^2\right)^{-\delta}\d x\right). 
$$
By Lemma \ref{L34} it is enough to prove that $S= (S(t))$ extends to a $C_0$-semigroup in the spaces  $\mathcal{L}^p_\delta$  and $L^p_\theta$  separately. The $C_0$-property of the semigroup $(S(t))$ in $\mathcal{L}^p_\delta$  can be shown using the method from  \cite{Peszat-Zabczyk1}. Therefore, it remains to prove the semigroup property in $L^p_\theta$.  We have 
\begin{align*}
\int_{\mathcal{O}}\left\vert S(t) \varphi (x)\right\vert ^p\rho^\theta(x) \d x&=\int_{\mathcal{O}}\left\vert \rho^{\frac{\theta+1}{p}}(x)S(t)\varphi(x)\right\vert ^p\frac{\d x}{\rho(x)}\\
&=\int_{\mathcal{O}}\left\vert \int_{\mathcal{O}}\rho^{\frac{\theta+1}{p}}(x)G(t, x,y)\rho(y)\varphi (y)\,\frac{\d y}{\rho(y)}\right\vert ^p\frac{\d x}{\rho(x)}\\
&=\int_{\mathcal{O}} \left\vert \int_{\mathcal{O}}\rho^{\frac{\theta+1}{p}}(x)G(t,x,y)\rho^{1-\frac{\theta+1}{p}}(y)\rho^{\frac{\theta+1}{p}}(y)\varphi (y)\,\frac{\d y}{\rho(y)}\right\vert ^p\frac{\d x}{\rho(x)}\\
&\le C \int_{\mathcal{O}}\left\vert \int_{\mathcal{O}}\rho^{\frac{\theta+1}{p}}(x) m_t(y)g_{ct}(x-y)\rho^{1- \frac{\theta+1}{p}}(y) \psi (y)\,\frac{\d y}{\rho(y)}\right\vert ^p \frac{\d x}{\rho(x)},
\end{align*}
with $\psi(y)=\rho^{\frac{\theta+1}{p}}(y)\varphi (y)$. Note that the last inequality follows from \eqref{E33}.  In other words 
$$
\left\vert S(t)\varphi \right\vert ^p_{L^p_\theta} \le  C \left\vert K_t \psi \right\vert ^p_{L^p\left({\mathcal{O}}, \frac{\d y}{\rho(y)}\right)}, 
$$
where 
\begin{align*}
K_t\psi(x)&:=\int_{{\mathcal{O}}}k_t(x,y)\psi(y)\,\frac{\d y}{\rho(y)}, \\
k_t(x,y)&:= \left( \frac{\rho(x)}{\rho(y)}\right)^{\frac{\theta+1}{p}}m_t(y)g_{ct}(x-y)\rho(y), 
\end{align*}
and $\psi$ is as above.  

Since $\varphi \mapsto \psi= \rho^{\frac{\theta+1}{p}}\varphi $ is an isometry between $L^p_\theta$ and $L^p\left({\mathcal{O}}, \frac{\d y}{\rho(y)}\right)$, the proof of a $C_0$-property will be completed as soon as we show that for each $0<t\le 1$, $K_t$ is a bounded linear operator from $L^p\left({\mathcal{O}}, \frac{\d y}{\rho(y)}\right)$ into $L^p\left({\mathcal{O}}, \frac{\d y}{\rho(y)}\right)$ and that $\sup_{0<t\le 1} \left\| K_t\right\| <+\infty$, where $\|\cdot \|$ is the operator norm on $L\left(L^p\left({\mathcal{O}}, \frac{\d y}{\rho(y)}\right), L^p\left({\mathcal{O}},  \frac{\d y}{\rho(y)}\right)\right)$.  The second part of the theorem follows since, by \eqref{E34},  
$$
\left\vert \frac{\partial S(t)\varphi}{\partial x_i}  \right\vert ^p_{L^p_\theta} \le  C t^{-\frac{p}{2}} \left\vert K_t \psi \right\vert ^p_{L^p\left({\mathcal{O}}, \frac{\d y}{\rho(y)}\right)}. 
$$

Taking into account the Schur test, see e.g. Theorem 5.9.2 in \cite{Garling},  it is enough to show that 
$$
\sup_{0<t\le 1}\sup_{x\in {\mathcal{O}}}\int_{\mathcal{O}}k_t(x,y)\,\frac{\d y}{\rho(y)}+ \sup_{0<t\le 1}\sup_{y\in{\mathcal{O}} }\int_{\mathcal{O}}k_t(x,y)\,\frac{\d x}{\rho(x)}<+ \infty\,.
$$
Note that  our  assumption $\theta <2p-1$ is necessary for the application of the Schur test. 

Given $t\in (0,1]$,  let ${\mathcal{O}}_t :=\left\{ x\in {\mathcal{O}}\colon \rho(x)< \sqrt{t}\right\}$ and $\left( {\mathcal{O}}_t\right)^{\mathrm{c}}:= \mathcal{O}\setminus {\mathcal{O}}_t$. Write 
\begin{align*}
k_1&:= \sup_{0<t\le 1}\sup_{x\in {\mathcal{O}}_t}\int_{{\mathcal{O}}_t}k_t(x,y)\,\frac{\d y}{\rho(y)}, &
k_2&:= \sup_{0<t\le 1}\sup_{y\in {\mathcal{O}}_t}\int_{{\mathcal{O}}_t}k_t(x,y)\,\frac{\d x}{\rho(x)},\\
k_3&:= \sup_{0<t\le 1}\sup_{x\in\left({\mathcal{O}}_t\right)^\mathrm{c}}\int_{{\mathcal{O}}_t}k_t(x,y)\,\frac{\d y}{\rho(y)},  & 
k_4&:= \sup_{0<t\le 1}\sup_{y\in {\mathcal{O}}_t }\int_{\left({\mathcal{O}}_t\right)^\mathrm{c}}k_t(x,y)\,\frac{\d x}{\rho(x)}, \\
k_5&:= \sup_{0<t\le 1}\sup_{x\in {\mathcal{O}}_t}\int_{\left({\mathcal{O}}_t\right)^\mathrm{c}}k_t(x,y)\,\frac{\d y}{\rho(y)}, &
k_6&:= \sup_{0<t\le 1}\sup_{y\in \left({\mathcal{O}}_t\right)^\mathrm{c}}\int_{{\mathcal{O}}_t}k_t(x,y)\,\frac{\d x}{\rho(x)}, \\
k_7&:= \sup_{0<t\le 1}\sup_{x\in\left({\mathcal{O}}_t\right)^\mathrm{c}}\int_{\left({\mathcal{O}}_t\right)^\mathrm{c}}k_t(x,y)\,\frac{\d y}{\rho(y)}, &
k_8&:= \sup_{0<t\le 1}\sup_{y\in\left({\mathcal{O}}_t\right)^\mathrm{c}}\int_{\left({\mathcal{O}}_t\right)^\mathrm{c}}k_t(x,y)\,\frac{\d x}{\rho(x)}.
\end{align*}
Note that the proof will be completed as soon as we show that all $k_j$ are finite. 

To estimate  $k_5$ to  $k_8$ where $y\in \left({\mathcal{O}}_t\right)^\mathrm{c}$ we use the  Lipchitz continuity of the distance function $\rho$ and the estimate $\rho(y)\ge \sqrt{t}$ for $y\in \left({\mathcal{O}}_t\right)^\mathrm{c}$. Namely,  for any $\alpha \ge 0$, we have 
$$
\left(\frac{\rho(x)}{\rho(y)}\right)^{\alpha} \le \left(\frac{\vert \rho(x)-\rho (y)\vert }{\rho(y)}+ 1 \right)^{\alpha} \le C \left(\frac{\vert x-y\vert }{\sqrt{t}}+ 1 \right)^{\alpha}.
$$
Since  $m_t(y)\le 1$, we have 
\begin{align*}
k_5&\le  \sup_{0<t\le 1}\sup_{x\in{\mathcal{O}}_t} \int_{\left({\mathcal{O}}_t\right)^\mathrm{c}} \left(\frac{\rho(x)}{\rho(y)}\right)^{\frac{\theta+1}{p}} g_{ct}(x-y) \d y\\
&\le  \sup_{0<t\le 1}\sup_{x\in {\mathcal{O}}}\int_{\left({\mathcal{O}}_t\right)^\mathrm{c}}\left(\frac{\vert x-y\vert }{\sqrt{t}}+ 1 \right)^{\frac{\theta +1}{p}} g_{ct}(x-y)\d y\\
&\le  \sup_{0<t\le 1}\int_{\mathbb{R}^d}\left(\frac{\vert z\vert }{\sqrt{t}}+ 1 \right)^{\frac{\theta +1}{p}} g_{ct}(z)\d z\le \int_{\mathbb{R}^d}\left(\vert z\vert + 1 \right)^{\frac{\theta +1}{p}-1} g_{c}(z)\d z<+\infty. 
\end{align*}
To estimate $k_7$ note that 
\begin{align*}
k_7 & \le  \sup_{0<t\le 1}\sup_{x\in\left({\mathcal{O}}_t\right)^\mathrm{c}} \int_{\left({\mathcal{O}}_t\right)^\mathrm{c}} \left(\frac{\rho(x)}{\rho(y)}\right)^{\frac{\theta +1}{p}} g_{ct}(x-y)\d y \le C \int_{\mathbb{R}^d}\left(\vert z\vert + 1 \right)^{\frac{\theta +1}{p}} g_{c}(z)\d z<+\infty. 
\end{align*}
In  the case of $\frac{\theta+1}{p}-1>0$, equivalently of   $\theta >p-1$,  one can use the same arguments to evaluate $k_6$ and $k_8$.  Namely, we have 
\begin{align*}
k_6&\le  \sup_{0<t\le 1}\sup_{y\in\left({\mathcal{O}}_t\right)^\mathrm{c}} \int_{{\mathcal{O}}_t} \left(\frac{\rho(x)}{\rho(y)}\right)^{\frac{\theta+1}{p}} g_{ct}(x-y)\left( \frac{\rho(x)}{\rho(y)}\right)^{-1} \d x\\
&\le  C \sup_{0<t\le 1} \sup_{y\in\left({\mathcal{O}}_t\right)^\mathrm{c}}  \int_{{\mathcal{O}}_t} \left( \frac{\vert x-y\vert }{\sqrt{t}} +1 \right)^{\frac{\theta+1}{p}-1}(x) g_{ct}(x-y)\d x<+\infty
\end{align*}
and 
\begin{align*}
k_8 & \le  \sup_{0<t\le 1}\sup_{y\in\left({\mathcal{O}}_t\right)^\mathrm{c}} \int_{\left({\mathcal{O}}_t\right)^{\mathrm{c}}} \left(\frac{\rho(x)}{\rho(y)}\right)^{\frac{\theta +1}{p}} g_{ct}(x-y)\left( \frac{\rho(x)}{\rho(y)}\right)^{-1} \d x\\
&\le C \sup_{0<t\le 1} \sup_{y\in\left({\mathcal{O}}_t\right)^\mathrm{c}}  \int_{\left({\mathcal{O}}_t\right)^\mathrm{c}} \left( \frac{\vert x-y\vert }{\sqrt{t}} +1 \right)^{\frac{\theta+1}{p}-1}(x) g_{ct}(x-y)\d x<+\infty 
\end{align*}

The case  of $\frac{\theta +1}{p}-1<0$ can be treated as follows 
\begin{align*}
k_6&\le  \sup_{0<t\le 1}\sup_{y\in\left({\mathcal{O}}_t\right)^\mathrm{c}} \int_{{\mathcal{O}}_t} \left(\frac{\rho(x)}{\rho(y)}\right)^{\frac{\theta+1}{p}} g_{ct}(x-y)\left( \frac{\rho(x)}{\rho(y)}\right)^{-1} \d x\\
&\le  \sup_{0<t\le 1}  \sup_{\sqrt{t}\le u\le 1} \sup_{y\colon \rho(y)=u} u^{1-\frac{\theta +1}{p}} \int_{{\mathcal{O}}_t} \rho^{\frac{\theta+1}{p}-1}(x) g_{ct}(x-y)\d x. 
\end{align*}
Note that  for any $\sqrt{t}\le u\le 1$, we have  
$$
\inf\left\{\vert x-y\vert ^2\colon y\in \left({\mathcal{O}}_t\right)^{\mathrm{c}}, \rho(y)=u, \ x\in {\mathcal{O}}_t \right\} =  \vert u-\sqrt{t}\vert^2.
$$
Thus, by Assumption \ref{Axx},   
\begin{align*}
k_6&\le \left(2\pi c\right)^{\frac d2} \sup_{0<t\le 1}\sup_{\sqrt{t}\le u\le 1} \sup_{y\colon \rho(y)=u} u^{1-\frac{\theta +1}{p}}\e^{-\frac{\vert u-\sqrt{t}\vert }{4ct}} \int_{{\mathcal{O}}_t} \rho^{\frac{\theta+1}{p}-1}(x) g_{2ct}(x-y)\d x\\
&\le \left(2\pi c\right)^{\frac d2} \sup_{0<t\le 1}\sup_{\sqrt{t}\le u\le 1} \sup_{y\colon \rho(y)=u}u^{1-\frac{\theta +1}{p}} \e^{-\frac{\vert u-\sqrt{t}\vert }{4ct}} t^{\frac{\theta+1}{2p}-\frac 12}<+\infty. 
\end{align*}
In the same way, if $\frac{\theta +1}{p}-1< 0$, then 
\begin{align*}
k_8 & \le  \sup_{0<t\le 1}\sup_{\sqrt{t}\le u\le 1} \sup_{y\colon \rho(y)=u} \int_{\left({\mathcal{O}}_t\right)^{\mathrm{c}}} \left(\frac{\rho(x)}{\rho(y)}\right)^{\frac{\theta +1}{p}} g_{ct}(x-y)\left( \frac{\rho(x)}{\rho(y)}\right)^{-1} \d x\\
&\le \left(2\pi c\right)^{\frac d2} \sup_{0<t\le 1}\sup_{\sqrt{t}\le u\le 1} \sup_{y\colon \rho(y)=u} u^{1-\frac{\theta +1}{p}}\e^{-\frac{\vert u-\sqrt{t}\vert }{4ct}} \int_{{\mathcal{O}}} \rho^{\frac{\theta+1}{p}-1}(x) g_{2ct}(x-y)\d x<+\infty. 
\end{align*}

We use Assumption \ref{Axx}  to evaluate $k_4$. Namely, since  $2-\frac{\theta +1}{p}>0$,   we have 
\begin{align*}
k_4&=  \sup_{0<t\le 1}\sup_{y\in {\mathcal{O}}_t} \rho^{2-\frac{\theta+1}{p}}(y)t^{-\frac12} \int_{ \left({\mathcal{O}}_t\right)^\mathrm{c}} \rho^{\frac{\theta +1}{p}-1}(x)  g_{ct}(x-y)\,\d x\\
&\le  \sup_{0<t\le 1} t^{\frac{1}{2}-\frac{\theta+1}{2p}} \sup_{y\in {\mathcal{O}}_t} \int_{ \left({\mathcal{O}}_t\right)^\mathrm{c}}\rho^{\frac{\theta +1}{p}-1}(x) g_{ct}(x-y)\,\d x<+\infty.
\end{align*}
The same argument can be used to evaluate $k_1$. Namely  since for $x\in {\mathcal{O}}_t$, $\rho(x)\le \sqrt{t}$, we have 
\begin{align*}
k_1 &= \sup_{0<t\le 1}\sup_{x\in {\mathcal{O}}_t}\int_{{\mathcal{O}}_t}\rho^{\frac{\theta+1}{p}}(x)\frac{\rho(y)}{\sqrt{t}}g_{ct}(x-y)\rho^{- \frac{\theta+1}{p}}(y)\,\d y\\
&\le C_1 \sup_{0<t\le 1}t^{\frac{\theta+1}{2p}-\frac{1}{2}}  \sup_{x\in {\mathcal{O}}_t}  \int_{{\mathcal{O}}_t}\rho^{1 - \frac{\theta+1}{p}}(y)g_{ct}(x-y)\d y<+\infty.
\end{align*}
Above we used  Assumption \ref{Axx} and the fact that $1-\frac{\theta +1}{p}>-1$ as $\theta < 2p-1$. 

To estimate $k_2$ we need  $2-\frac{\theta +1}{p}>0$, that is $\theta < 2p-1$.  Since  $\rho(x)\le \sqrt{t}$ and $\rho(y)\le \sqrt{t}$  for $x,y\in {\mathcal{O}}_t$, we have 
\begin{align*}
k_2 &= \sup_{0<t\le 1}\sup_{y\in {\mathcal{O}}_t}\int_{{\mathcal{O}}_t}\rho^{\frac{\theta+1}{p}}(x) t^{-1}g_{ct}(x-y)\rho^{2- \frac{\theta+1}{p}}(y)\,\d x\\
&\le  \sup_{0<t\le 1}  t^{\frac{\theta +1}{2p}-1 + 1-\frac{\theta +1}{2p}}  \sup_{y\in {\mathcal{O}}_t} \int_{{\mathcal{O}}_t}g_{ct}(x-y) \d x<+\infty. 
\end{align*}

It remains to evaluate  $k_3$. We have 
\begin{align*}
k_3&= \sup_{0<t\le 1}\sup_{x\in \left({\mathcal{O}}_t\right)^\mathrm{c}} \rho^{\frac{\theta +1}{p}}(x)t^{-\frac{1}{2}} \int_{\mathbb{B}^d_t} g_{ct}(x-y)\rho^{1-\frac{\theta +1}{p}}(y)\,\d y. 
\end{align*}
Note that  for any $\sqrt{t}\le u\le 1$, we have  
$$
\inf\left\{\vert x-y\vert ^2\colon x\in \left({\mathcal{O}}_t\right)^{\mathrm{c}}, \rho(x)=u, \ y\in {\mathcal{O}}_t \right\} =  \vert u-\sqrt{t}\vert^2.
$$
Thus, by Assumption \ref{Axx},   
\begin{align*}
k_3&\le C_1 \sup_{0<t\le 1}\sup_{\sqrt{t}0\le u\le 1} \sup_{x\in {\mathcal{O}}_t\colon \rho(x)=u} u ^{\frac{\theta +1}{p}}t^{-\frac{1}{2}}\e^{-\frac{\vert u-\sqrt{t}\vert^2}{4ct}} \int_{{\mathcal{O}}_t} \rho^{1-\frac{\theta +1}{p}}(y)g_{2ct}(x-y)\d y\\
&\le C_2 \sup_{0<t\le 1}\sup_{\sqrt{0}\le u\le 1} u ^{\frac{\theta +1}{p}}t^{-\frac{1}{2}}\e^{-\frac{\vert u-\sqrt{t}\vert^2}{4ct}} t^{\frac 12 - \frac{\theta+1}{2p}}<+\infty. 
\end{align*}
$\square$
\subsection{Analiticity}
\begin{remark}
{\rm Assume that the derivatives $\frac{\partial }{\partial x_i}$ commute with the semigroup in the following sense 
$$
\frac{\partial }{\partial x_i}S(t)= S(t/2)\frac{\partial }{\partial x_i}S(t/2) + R_i(t)S(t/2),
$$
where $R_i(t)$, $t>0$ are  bounded linear operator satisfying 
$$
\|R_i(t)\|_{L(L^p_{\theta,\delta},L^p_{\theta,\delta})}\le C_1t^{-1/2}. 
$$
Then, by second part of Theorem \ref{T34}, 
\begin{align*}
\left\vert \frac{\partial ^2}{\partial x_i^2} S(t)\psi\right\vert _{L^p_{\theta,\delta}}&= \left\vert \frac{\partial }{\partial x_i} S(t/2)\frac{\partial }{\partial x_i}S(t/2)\psi+ \frac{\partial }{\partial x_i} S(t/2) R_i(t)S(t/2) \psi \right\vert _{L^p_{\theta,\delta}}\le \frac{C^2}{t} \left\vert \psi\right\vert _{L^p_{\theta,\delta}}. 
\end{align*}
This leads to the analiticity of $S$ on $L^p_{\theta,\delta}$ in the case of $\mathcal{A}$ of the form $\sum_{i,j} a_{i,j} \frac{\partial ^2}{\partial x_i\partial x_j}+ \sum_i b_i \frac{\partial }{\partial x_i}$ with bounded $a_{i,j}$ and $b_j$.  }
\end{remark}

The classical Aronson estimates for the Green kernel,  see e.g. \cite{Eidelman-Ivasishen, Solonnikov, Mora} for required assumptions on $\mathcal{A}$ and $\mathcal{O}$, yield  that  $G$ is of class $C^\infty ((0,+\infty)\times \mathcal{O}\times \mathcal{O})$ and for any non-negative integer $n$,  multi-indices  $\alpha $, $\beta$,  and time $T>0$, there are constants $C,c>0$ such that for all $t\in (0,T]$ and $x, y\in \mathcal{O}$,
$$
\left \vert \frac{\partial^n}{\partial t^n} \frac{\partial ^{|\alpha|}}{\partial x^\alpha }\frac{\partial ^{|\beta|}}{\partial y^\beta }  G(t,x,y)\right\vert \le C  t ^{-\frac{|\alpha|+|\beta| + 2n}{2}} g_{ct}(x-y).
$$
In our  proofs of the $C_0$-property and  gradient estimate we needed something different, namely estimates  \eqref{E33} and \eqref{E34} which guarantee that $G(t,x,y)$ and $\nabla _xG(t,x,y)$ decay for $y$ near the boundary of $\mathcal{O}$ at rate $\rho(y)/\sqrt{t}$ uniformly in $x$. Clearly, our proof yelds the following. 
\begin{proposition}  If for a certain multi index $\alpha$, there are constants $C,c>0$ such that   
\begin{equation}\label{Ede}
\left \vert\frac{\partial ^{|\alpha|}G}{\partial x^\alpha }(t,x,y)\right\vert \le C t^{-\frac{\vert \alpha \vert}{2}}m_t(y)g_{ct}(x-y), \qquad \forall\, x,y\in \mathcal{O},\ \forall\, t\in (0,1],
\end{equation} 
then, for all $p\ge 1$, $\theta \in [0, 2p-1)$, $\delta \ge 0$, and $T>0$, there is a constant $C_1$ such that 
$$
\left\vert \frac{\partial ^{\vert \alpha \vert}}{\partial x^\alpha } S(t)\psi\right \vert _{L^p_{\theta,\delta}}\le C_1 t^{-\frac{\vert \alpha \vert}{2}}\left\vert \psi\right\vert _{L^p_{\theta,\delta}},\qquad \forall\, \psi\in L^p_{\theta,\delta},\ t\in (0,T]. 
$$
\end{proposition}
\begin{corollary}
If $\mathcal {A}=\Delta$ and $\mathcal{O}$ is a half space, then for all $1\le p<+\infty$, $\theta \in [0,2p-1)$ and $\delta\ge 0$, the semigroup $S$ is analytical on  $L^p_{\theta,\delta}$. 
\end{corollary}
\begin{proof} We need to show that there is a constant $C$ such that 
$$
\left\vert \Delta S(t)\psi\right\vert _{L^p_{\theta,\delta}}\le \frac{C}{t}\left\vert \psi\right\vert _{L^p_{\theta,\delta}} ,\qquad \forall\, \psi\in L^p_{\theta,\delta},\ t\in (0,1]. 
$$
We may assume that $\mathcal{O}=\{x\in \mathbb{R}^d\colon x_1>0\}$. Then $m_t(y)= (y_1/\sqrt{t})\wedge 1$, and  the Green kernel  in known, namely 
\begin{equation}\label{EK1}
G(t,x,y)= g_{2t} (x-y)-g_{2t} (\overline x-y),
\end{equation}
 where 
\begin{equation}\label{EK2}
\overline {x}=\overline{(x_1,x_2,\ldots,x_d)}=\overline {(x_1, \mathbf{x})} = (-x_0,\mathbf{x}).
\end{equation}
By elementary calculation one can verified  estimate \eqref{Ede} for any second order derivative $\frac{\partial ^2}{\partial x_j^2}$. Indeed, given  $a>0$, $z\in \mathbb{R}$ and $\mathbf{z}\in \mathbb{R}^{d-1}$ write 
$$
g^1_{a}(z):= \left(2\pi a\right)^{-\frac 12} \e ^{-\frac {z^2}{2a}}, \qquad g^{d-1}_{a}(\mathbf{z}):= \left(2\pi a\right)^{-\frac {d-1}2} \e ^{-\frac {\vert \mathbf{z}\vert ^2}{2a}}. 
$$
Note that  there is a constant $C$ such that for all $x_1,y_1\ge 0$, $t\in (0,1]$, 
\begin{equation}\label{Etr}
\left\vert  g^1_{2t}(x_1-y_1)- g^1_{2t}(x_1+y_1)\right\vert \le C m_t(y) g^1_{4t}(x_1-y_1). 
\end{equation}
For, \eqref{Etr} can be reformulated equivalently as 
$$
\left\vert \e^{-z^2}-\e^{-(z+v)^2}\right\vert\le C\, v\wedge 1\, \e^{-\frac{z^2}{2}}, \qquad \forall\, z\in \mathbb{R}, v \ge 0. 
$$
or 
$$
\e^{-\frac {z^2}{2}} \left\vert 1-\e^{-(z+v)^2+z^2 }\right\vert\le C\, v\wedge 1, \qquad \forall\, z\in \mathbb{R}, v \ge 0. 
$$
We have 
\begin{align*}
\left\vert \frac{\partial^{2} G}{\partial x_1^2}(t,x,y)\right\vert &= \left\vert \left[ \frac{(x_1-y_1)^2}{4t^2}-\frac 1{2t}\right] g_{2t}(x-y) +  \left[- \frac{(x_1+y_1)^2}{4t^2}+\frac 1{2t}\right] g_{2t}(\overline x-y)\right\vert. 
\end{align*}
Therefore, by \eqref{Etr}, it is enough to show that for all $u,v\ge 0$, 
\begin{align*}
\left\vert (u-v)^2 \e^{-\frac{(u-v)^2}{4}}- (u+v)^2\e^{-\frac{(u+v)^2}{4}}\right\vert&\le  C\, v\wedge 1 \, \e^{-\frac{(u-v)^2}{8}}. 
\end{align*}

For $j>1$ we have 
\begin{align*}
\left\vert \frac{\partial^{2} G}{\partial x_j^2}(t,x,y)\right\vert &= \left\vert  g^1_{2t}(x_1-y_1)- g^1_{2t}(x_1+y_1)\right\vert  \left\vert  \frac{(x_j-y_j)^2}{4t^2 }- \frac {1}{2t}\right\vert g^{d-1}_{2t}(\mathbf{x}-\mathbf{y}) \\
&\le \frac{C_1}{t} \left\vert  g^1_{2t}(x_1-y_1)- g^1_{2t}(x_1+y_1)\right\vert  g^{d-1}_{4t}(\mathbf{x}-\mathbf{y})\\
&\le  \frac{CC_1}{t} m_t(y)   g_{4t}(x-y). 
\end{align*}
\end{proof}
 \subsection{Related  results}
 In this section we comment some recent results of Krylov \cite{Krylov} and \cite{Krylov} and Lindemulder and Veraar \cite{Lindemulder-Veraar} concerning heat semigroup on weighted spaces. 
 \subsubsection{Krylov's result}
Let $P=\left\{ x\in \mathbb{R}^d\colon x_1> 0\right\}$ be a half space in $\mathbb{R}^d$. Let $\rho(x)= x_1$ be the distance of $x\in P$ from the boundary. Let $S$ be the semigroup generated by the Laplace operator $A=\Delta$ on $P$ with homogeneous Dirichlet boundary conditions. By $\nabla $ we denote the gradient operator and by $\nabla ^2$ the Hessian.  Given $\theta \in \mathbb{R}$ let $L^p_\theta= L^p(P,x_1^\theta \d x)$.  The following result follows directly from the Krylov Theorem 2.5  (\cite{Krylov}). In the original Krylov theorem $p=q$, $\alpha =2=\hat\alpha =a_+$, $\gamma = \overline \gamma =0$. 
\begin{theorem}
Let $p\in (1,+\infty)$. Then for every $\theta \in (-2p,p)$, $S$ is a $C_0$-semigroup on $L^p_\theta$. Moreover, there is a constant $N$ such that   for any $t>0$, 
$$
\|S(t)u \|_{L^p_\theta}\le N \|u\|_{L^p_{\theta}}\qquad \text{and}\qquad   \|\nabla ^2 S(t)u \|_{L^p_\theta}\le Nt^{-1}\|u\|_{L^p_{\theta}}. 
$$
\end{theorem}
Given a vector $a\in \mathbb{R}^d$ and a number $\delta \in \mathbb{R}$ let us denote by $P(a,\delta)$ the half space 
$$
P(a,\delta ):=\left\{x\in \mathbb{R}^d\colon \langle x,a\rangle >\delta\right\}.
$$
Let $\rho_{P(a,\delta)}(x)$ be the distances of $x\in \mathbb{R}^d$ from the boundary $\partial P(a,\delta)$.  Obviously the Krylov result can be extended to  any of half space  $P(a,\delta)$. The $L^p_\theta$ space should be replaced by 
$$
L^p_\theta(P(a,\delta )):=L^p(P(a,\delta),\rho^\theta _{P(a,\delta)}(x)\d x). 
$$ 
Note that the constant $N$ appearing in the theorem is universal for any half space.

Let $\mathcal{ O}$ be a not necessarily bounded domain in $\mathbb{R}^d$. Let $\rho_{\mathcal{O}}(x)$ be the distance of $x\in \mathcal{ O}$ from the boundary $\partial \mathcal{ O}$. Given $\theta \in \mathbb{R}$ write $L^p_{\theta}= L^p(\mathcal{ O}, \rho^\theta_{\mathcal{ O}}(x)\d x)$. Let $S$ be the heat semigroup on $\mathcal{ O}$ with homogeneous Dirichlet boundary conditions. 
\begin{theorem}
Assume that $\mathcal  O$ is  a convex domain in $\mathbb{R}^d$. Let $p\in (1,+\infty)$. Then for every $\theta \in (-2p,p)$, $S$ is a $C_0$-semigroup on $L^p_\theta$. Moreover, there is an independent of $\mathcal  O$  constant $N$ such that   for any $t>0$, 
\begin{equation}\label{E39}
\vert S(t)\psi \vert _{L^p_\theta}\le N \vert \psi\vert _{L^p_{\theta}}, \qquad \psi\in L^p_\theta. 
\end{equation}
\end{theorem}
\begin{proof} Since $\mathcal  O$ is convex then there is a family of subspaces $P(a_j,\delta_j)$, $j\in J$, such that 
\begin{equation}\label{E310}
\mathcal O=\bigcap_{j\in J} P(a_j,\delta _j). 
\end{equation}
Let $j\in J$, and let $\psi\in C^\infty_{0}(\mathcal  O)$. Let $T_{a_j,\delta_j}$ be the heat semigroup on $P(a_j,\delta_j)$ with homogeneous Dirichlet boundary conditions.  Let us observe that 
\begin{equation}\label{E311}
\vert S(t) \psi (x)\vert \le T_{a_j,\delta_j }(t)\vert \psi \vert (x), \qquad x\in \mathcal O. 
\end{equation}
For, \eqref{E311}  follows immediately for example from the following probabilistic representations 
\begin{align*}
S(t)\psi (x)&= \mathbb{E}\left(\psi (x+W(t)); t<\tau_x(\mathcal O)\right),\\ T_{a_j,\delta_j}(t)\psi (x)&= \mathbb{E}\left(\psi (x+W(t)); t<\tau _x\left(P(a_j,\delta_j)\right)\right),
\end{align*}
where $\tau_x(\mathcal  O)$ and $\tau _x\left(P(a_j,\delta_j)\right)$ are exit times 
$$
\tau_x(\mathcal  O):= \inf\{s>0\colon x+W(s)\not \in \mathcal  O\}, \qquad \tau _x\left(P(a_j,\delta_j)\right):=\inf\{s>0\colon x+W(s)\not \in P(a_j,\delta_j)\}. 
$$
Thus 
$$
\left\vert   S(t) \psi (x)\right\vert ^p \rho^\theta _{\mathcal  O}(x)\le \left\vert   T_{a_j,\delta_j}(t) \psi (x)\right\vert ^p \rho^\theta _{P(a_j,\delta_j)}(x), \qquad \forall\, x\in \mathcal  O. 
$$
Therefore after integration we obtain 
$$
\vert S(t)\psi \vert _{L^p_\theta}\le \vert T_{a_j,\delta_j}(t)\psi \vert _{L^p_\theta(P(a_j,\delta_j))}\le N \vert \psi \vert _{L^p_\theta(P(a_j,\delta_j))}, 
$$
where the constant $N$ does not depend on $j$. Hence 
$$
\vert   S(t) \psi \vert _{L^p_ \theta}\le N \inf_{j\in J} \vert \psi \vert _{L^p_\theta(P(a_j,\delta_j))}. 
$$
Taking into account \eqref{E310} we have 
$$
\inf_{j\in J} \vert \psi \vert _{L^p_\theta(P(a_j,\delta_j))}= \vert \psi \vert _{L^p_\theta}, 
$$ 
which gives \eqref{E39} and obviously $C_0$-property of $S$.  
\end{proof}
\begin{remark}
{\rm Unfortunately, since we do not have the estimate for the gradient 
$$
\vert \nabla S(t) \psi (x)\vert \le \vert \nabla T_{a_j,\delta_j }(t) \psi  (x)\vert , \qquad x\in \mathcal  O.
$$
the derivation of  the estimate  $\vert \nabla  S(t)\psi  \vert _{L^p_\theta}\le Nt^{-1/2}\vert \psi\vert _{L^p_{\theta}}$ needs some different arguments! }
\end{remark}
\subsubsection{Lindemulder and Veraar results} As in the Krylov papers, paper \cite{Lindemulder-Veraar} of Lindemulder and Veraar deals with the Laplace operator $\mathcal{A}=\Delta$. It is shown that $\Delta$ with Dirichlet boundary conditions admits  a bounded $H^\infty$-calculus on weighted spaces $L^p_\theta:=L^p\left(\mathcal{O},\rho^\theta(x)\d x\right)$, where $\rho(x)=  \textrm{dist}\left(x,\partial \mathcal{O}\right)$, $p\in (1,+\infty)$ and $\theta \in (-1,2p-1)\setminus \{p-1\}$.  Therefore, the corresponding heat semigroup is not only $C_0$ but also analytical on  $L^p_\theta$.   In \cite{Lindemulder-Veraar}, $\mathcal{O}$ is a halfspace or a bounded $C^2$-domain.  

\section{Properties of the semigroup on weighted spaces}\label{SStability}
In this section, Assumptions \ref{A31} and \ref{A32} are satisfied. 
 
 \begin{lemma}\label{L41}
 There exists a $C>0$ such that 
 $$
 \left|S(t)\psi\right |_{\mathcal{L}^p_\delta}\le Ct^{-\frac{\theta}{2p}}|\psi|_{L^p_{\theta,\delta}}  \qquad  \text{for $t\in(0,1]$.} 
$$
 \end{lemma}
 \begin{proof}  Let
 \begin{equation}\label{E41}
 w_\delta(x)= w_{0,\delta}(x)= \left(1 + |x|^2\right)^{-\delta}
 \end{equation}
 be the weight on $\mathcal{L}^p_\delta$. Then, by Assumption \ref{A32},  we have 
$$
\begin{aligned}
\left | S(t)\psi \right|^p _{\mathcal{L}^p_\delta}&\le C^p  \int_{\mathcal O} w_\delta(x) \left(  \int_{\mathcal{O}}m_t(y)g_{ct}(x-y)\left|\psi(y)\right| \d y\right )^p \d x\\
&\le  \tilde C t^{d/2}(I_1+I_2),
\end{aligned}
$$
where 
$$
\begin{aligned}
I_1&:= \int_{\mathcal {O}/\sqrt{t}}w_\delta(x\sqrt{t})\left( \int_{\left( \mathcal{O}/\sqrt{t}\right)_1}m_t(y\sqrt{t})g_{c}(x-y)\left| \phi(y)\right| \d y\right)^p \d x, \\
I_2&:= \int_{\mathcal{O}/\sqrt{t}}w_\delta(x\sqrt{t})\left( \int_{\left(\mathcal{O}/\sqrt{t}\right)_1^c}m_t(y\sqrt{t})g_{c}(x-y)\left| \phi(y)\right| \d y\right)^p \d x,
 \end{aligned}
$$
$\phi(z)=\psi(z\sqrt{t})$, and 
$$
 (\mathcal{O}/\sqrt{t})_1:= \left\{x\in \mathcal{O}/\sqrt{t} \colon m_t(x\sqrt{t})=\rho(x\sqrt{t})/\sqrt{t}<1\right\}, 
$$
and 
$$
(\mathcal{O}/\sqrt{t})_1^c:= \left\{x\in \mathcal{O}/\sqrt{t}\colon m_t(x\sqrt{t})= 1\right\}. 
$$

We have 
$$
\begin{aligned}
I_1&\le  \int_{\mathcal{O}/\sqrt{t}}w_\delta(x\sqrt{t})\int_{\left( \mathcal{O}/\sqrt{t}\right)_1}m_t^p(y\sqrt{t})g_c(x-y)|\phi(y)|^p\d y\d x\\
&\le \int_{\left( \mathcal{O}/\sqrt{t}\right)_1}\left(\int_{\mathcal{O}/\sqrt{t}}w_\delta(x\sqrt{t})\dfrac{m_t^p(y\sqrt{t})}{w_{\theta,\delta}(y\sqrt{t})}g_c(x-y)\d x\right) \left |\phi(y)\right |^pw_{\theta,\delta}(y\sqrt{t})\d y \\
&\le \int_{\left( \mathcal{O}/\sqrt{t}\right)_1}F(t,y)\left |\phi(y)\right |^pw_{\theta,\delta}(y\sqrt{t})\d y,
\end{aligned}
$$
where 
\begin{align*}
F(t,y)& =\int_{\mathcal{O}/\sqrt{t}}w_\delta(x\sqrt{t})\dfrac{m_t^p(y\sqrt{t})}{w_{\theta,\delta}(y\sqrt{t})}g_c(x-y)\d x,\\
&= \int_{\mathcal{O}/\sqrt{t}}w_\delta(x\sqrt{t})\dfrac{\rho^p(y\sqrt{t})}{t^{p/2}w_{\theta,\delta}(y\sqrt{t})}g_c(x-y)\d x, \qquad y\in\left(\mathcal{ O}/\sqrt{t}\right)_1\,.
\end{align*}
Recall that $w_{\theta,\delta}(z)= \min\left\{ \rho^\theta(z), w_{\delta}(z)\right\}$. Thus, if  $\rho^\theta(y\sqrt{t})\le w_\delta(y\sqrt{t})$ then since 
$$
\rho^{p-\theta}(y\sqrt{t})\le t^{(p-\theta)/2}\quad \text{for}\quad y\in\left(\mathcal{ O}/\sqrt{t}\right)_1,  
$$
we have 
$$
\begin{aligned}
F(t,y)&\le \int_{\mathcal{O}/\sqrt{t}}\dfrac{\rho^{p-\theta}(y\sqrt{t})}{t^{p/2}}g_c(x-y) w_\delta(x\sqrt{t})\d x \le\dfrac{C_1}{t^{\theta/2}}\,.
\end{aligned}
$$
If $ w_\delta(y\sqrt{t})<\rho^\theta(y\sqrt{t})$ then 
$$
\begin{aligned}
F(t,y)&\le\int_{ \mathcal{ O}/\sqrt{t}}\dfrac{w_\delta(x\sqrt{t})}{w_{\delta}(y\sqrt{t})}g_c(x-y)\d x\,.
\end{aligned}
$$
Putting $a=\sqrt{t}$ we find that 
$$
\begin{aligned}
\dfrac{1+a|y|^2}{1+a|x|^2}&\le\dfrac{1+2a|y-x|^2 +2a|x|^2}{1+a|x|^2}\\
&\le 1+2a|x-y|^2+\dfrac{2a|x|^2}{1+a|x|^2}\\
&\le 3+2a|x-y|^2. 
\end{aligned}
$$
Therefore 
$$
\sup_{y\in\left( \mathcal{ O}/\sqrt{t}\right)_1}F(t,y)\le \frac{C_2}{t^{\theta/2}}\,,
$$
and hence
$$
t^{\frac d2}I_1\le C_2 t^{-\theta/2}|\psi |^p_{L^p_{\theta,\delta}}.
$$
For $I_2$ we obtain 
$$
\begin{aligned}
I_2&=\int_{\mathcal{ O}/\sqrt{t}}w_\delta(x\sqrt{t})\left(\int_{\left( \mathcal{ O}/\sqrt{t}\right)_1^c}m_t(y\sqrt{t})g_c(x-y)\phi(y)\d y\right) ^p\d x\\
&\le \int_{\mathcal{ O}/\sqrt{t}}w_\delta(x\sqrt{t})\int_{\left( \mathcal{ O}/\sqrt{t}\right)_1^c}\dfrac{g_c(x-y)}{w_{\theta,\delta}(y\sqrt{t})}|\phi(y)|^pw_{\theta,\delta}(y\sqrt{t})\d y\d x\\
&\le \int_{\left(\mathcal{ O}/\sqrt{t}\right)_1^c} H(t,y)|\phi(y)|^pw_{\theta,\delta}(y\sqrt{t})\d y, 
\end{aligned}
$$
where 
$$
H(t,y):=  \int_{\mathcal{ O}/\sqrt{t}}w_\delta(x\sqrt{t})\dfrac{g_c(x-y)}{w_{\theta,\delta}(y\sqrt{t})}\d x, \qquad y\in \left(\mathcal{ O}/\sqrt{t}\right)_1^c.
$$
Note that 
$$
\frac{\rho^\theta(y\sqrt{t})}{t^{\theta/2}}\ge 1. 
$$
Thus 
$$
1\ge w_{\theta,\delta}(y\sqrt{t})= \min\left\{ \rho^\theta(y\sqrt{t}),w_\delta(y\sqrt{t})\right\}\ge \min\left\{ t^{\theta/2},w_\delta(y\sqrt{t})\right\}. 
$$
Therefore 
\begin{align*}
H(t,y) &\le \int_{\mathcal{ O}/\sqrt{t}}\left[ \dfrac{w_\delta(x\sqrt{t})}{w_{\delta}(y\sqrt{t})} + w_\delta(x\sqrt{t})t^{-\theta/2}\right] g_c(x-y) \d x\\
&\le  \int_{\mathcal{ O}/\sqrt{t}}\left[ \dfrac{w_\delta(x\sqrt{t})}{w_{\delta}(y\sqrt{t})} + t^{-\theta/2}\right] g_c(x-y) \d x\\
&\le Ct^{-\theta /2}. 
\end{align*}
 \end{proof}
 In what follows we denote by $L(E,V)$ the space of all linear bounded  operators from a Banach  space $E$ to a Banach space $V$, equipped with the operator norm $\|\cdot \|_{L(E,V)}$. 
 \begin{theorem}\label{T42}
 Assume that there are positive constants $M,\alpha$ such that for every $t>0$, 
 $$
 \|S(t)\|_{L(\mathcal{L}^p_{\delta}, \mathcal{L}^p_{\delta})}\le M\e^{-\alpha t}\,.
 $$
Then there exist $C>0$, such that for every $t>0$,
$$
\|S(t)\|_{L(L^p_{\theta,\delta},  L^p_{\theta,\delta})}\le C\e^{-\alpha t}\,.
$$
 \end{theorem}
 \begin{proof} Clearly 
 $$
 \|S(t)\psi\|_{L^p_{\theta,\delta}} \le \|S(t)\psi\|_{\mathcal{L}^p_{\delta }}=  \|S(t-1)S(1)\psi\|_{\mathcal{L}^p_{\delta}}\le M\e^{-\alpha (t-1)} \|S(1)\psi\|_{\mathcal{L}^p_{\delta}}.
 $$
 By Lemma \ref{L41}, 
 $$
 \|S(1)\psi\|_{\mathcal{L}^p_{\delta}}\le C_2\|\psi\|_{L^p_{\theta,\delta}}. 
 $$
 \end{proof}
\section{Dirichlet map}\label{SDirichlet}
In this section we derive useful estimates for 
$$
\left( \lambda-A\right)S(t)D_\lambda e = S(t)Be, \qquad t\ge 0, e\in L^2(\partial {\mathcal O},\d \mathrm{s}). 
$$
We will need Assumptions \ref{A31} and \ref{A32}, and additionally the following assumption, which is satisfied, see Remark \ref{R33} if the drift coefficients $\nu^i$ of $\mathcal{A}$ are of the class $C^1_b$.
\begin{assumption}\label{Ass61}
Assume that the coefficients $a_{i,j}$ are bounded, and that for  any $T>0$ there is a constant $C$ such that 
$$
\left|\nabla_yG(t,x,y)\right|\le \frac{C}{{\sqrt{t}}}g_{ct}(x-y), \qquad t\le T, \qquad x,y\in \mathcal{O}. 
$$
\end{assumption}

Let $\mathcal{S}_0(\overline{\mathcal{O}})$ be the set of all $\psi$ of tempered test functions such that $\psi(x)=0$ for $x\in \partial \mathcal{O}$. Let $\gamma $ be a continuous compactly supported function on $\partial \mathcal{O}$. Recall   that  $u=D_\lambda \gamma$ is the solution to the non-homogeneous Poisson problem  \eqref{E13}.  Let $\psi\in \mathcal{S}_0(\overline{\mathcal{O}})$.  Then applying Gauss--Green integration by parts formula we obtain 
$$
\int_{\mathcal{O}}\mathcal{A}u(x)\psi(x)\d x = \int_{\mathcal{O}}u(x)\mathcal{A}^*\psi(x)\d x - \int_{\partial \mathcal{O}}\gamma(x)\sum_{i,j}a_{ij}(x)\frac{\partial \psi}{\partial x_i}(x)\mathbf{n}^j(x)\d \mathrm{s}(x),
$$
where $\mathcal{A}^*$ is the formal adjoint operator; 
$$
\mathcal{A}^*\psi(x)= \sum_{i,j}\frac{\partial }{\partial x_j} \left( a_{i,j}(x)\frac{\partial \psi}{\partial x_i}(x)\right)- \sum_{i} \frac{\partial }{\partial x_i}\left( \mu^i(x)\psi(x)\right),
$$
$\mathbf{n}=\left(\mathbf{n}^1,\ldots,\mathbf{n}^d\right)$ is the outward pointing unit normal vector to the boundary $\partial \mathcal{O}$,  and $\mathrm{s}$ is the surface measure.   

Therefore we have 
\begin{equation}\label{E52}
 \int_{\mathcal{O}}u(x)\mathcal{A}^*\psi(x)\d x = \lambda \int_{\mathcal{O}} u(x)\psi(x)\d x + \int_{\partial \mathcal{O}}\gamma(x)\sum_{i,j}a_{ij}(x)\frac{\partial \psi}{\partial x_i}(x)\mathbf{n}^j(x)\d\mathrm{s} (x). 
 \end{equation}
In fact \eqref{E52} can be treated  as the  definition of the weak solution to \eqref{E14}, see e.g. \cite{bgpr, Lions-Magenes}. 

Let 
$$
\mathbf{n}^a= \left(\sum_{j} a_{1,j}\mathbf{n}^j, \ldots, \sum_{j} a_{d,j}\mathbf{n}^j\right). 
$$
Then 
$$
\sum_{i,j}a_{ij}(x)\frac{\partial \psi}{\partial x_i}(x)\mathbf{n}_j(x)= \frac{\partial \psi}{\partial \mathbf{n}^a}(x),
$$
and \eqref{E52} has the form 
$$
 \int_{\mathcal{O}}D_\lambda\gamma(x)\left( \mathcal{A}^*-\lambda\right) \psi(x)\d x =  \int_{\partial \mathcal{O}}\gamma(x)\frac{\partial \psi}{\partial \mathbf{n}^a}(x)\d \mathrm{s}(x). 
$$

In what follows $A$ is the generator of the semigroup $S$ on $L^2(\mathcal{O})$,  $S^*$ is the adjoint semigroup and $A^*$ is its generator.  Note that $\mathcal{A}^*\subset A^*$. 
 \begin{proposition}\label{P51}
Let $\lambda$ be in the resolvent of $A^*$. Then the Dirichlet map is uniquely characterized by the relation 
 \begin{equation}\label{E53}
 \int_{\mathcal{O}} D_\lambda \gamma (x)\psi(x)\d x=  \int_{\partial \mathcal{O}}\gamma(x)\frac{\partial }{\partial \mathbf{n}^a}\left(A^*-\lambda\right)^{-1}\psi(x)\d \mathrm{s}(x), 
 \qquad   \psi\in \mathcal{S}_0(\overline{\mathcal{O}}). 
\end{equation}
\end{proposition}
\begin{proof} Assume \eqref{E53}. Let $\psi\in \mathcal{S}_0(\overline{\mathcal{O}})$. Then
\begin{align*}
\int_{\mathcal{O}}D_\lambda \gamma(x)\left( {\mathcal{A}}^*-\lambda\right)\psi(x)\d x &=\int_{\partial \mathcal{O}}\gamma(x) \frac{\partial }{\partial \mathbf{n}^a}\left(A^*-\lambda \right)^{-1}\left( \mathcal{A}^*-\lambda\right)\psi(x)\d \mathrm{s}(x)\\
&=  \int_{\partial \mathcal{O}}\gamma(x) \frac{\partial }{\partial {\mathbf{n}^a}} \psi(x)\d \mathrm{s}(x). 
\end{align*}
\end{proof}
 \begin{corollary}\label{C52}
We have 
\begin{equation}\label{E54}
\int_{\mathcal{O}}\left(\lambda -A\right) S(t) D_\lambda \gamma (x)\psi(x)\d x= - \int_{\partial \mathcal{O}}\gamma(x) \frac{\partial }{\partial \mathbf{n}^a}\left[ S^*(t)\psi(x)\right] \d \mathrm{s}(x). 
\end{equation}
\end{corollary}

Let $G$ be the Green kernel corresponding to the heat semigroup $S$ generated by $\mathcal{A}$ with homogeneous boundary conditions.  Let 
\begin{equation}\label{E55}
\mathcal{G}_\lambda (x,y)= \int_0^{+\infty} \e ^{-\lambda t} G(t,x,y)\d t, 
\end{equation}
where $\lambda$ is from the resolvent set. 
\begin{theorem}\label{T53}
We have: 
\begin{equation}\label{E56}
D_\lambda \gamma(x)= \int_{\partial \mathcal{O}}\gamma(y)\frac{\partial }{\partial \mathbf{n}^a(y)}\mathcal{G}_\lambda(x,y) \d \mathrm{s}(y)
\end{equation} 
and 
\begin{equation}\label{E57}
(\lambda -A)S(t)D_\lambda \gamma(x)= -\int_{\partial \mathcal{O}} \gamma(y) \frac{\partial }{\partial {\mathbf{n}^a}(y)} G(t,x,y)\d \mathrm{s}(y). 
\end{equation}
\end{theorem}
\begin{proof} Let $\psi\in \mathcal{S}_0(\overline{\mathcal{O}})$. Then, by \eqref{E52}, 
\begin{align*}
\int_{\mathcal{O}} D_\lambda \gamma (x)\psi(x)\d x&=  \int_{\partial \mathcal{O}}\gamma(x)\frac{\partial }{\partial \mathbf{n}^a}\left(A^*-\lambda\right)^{-1}\psi(x)\d \mathrm{s}(x)\\
 &= \int_{\partial \mathcal{O}}\gamma(x)\frac{\partial }{\partial \mathbf{n}^a(x)}\int_{\mathcal{O}}\mathcal{G}_\lambda(y,x)\psi(y)\d y \d \mathrm{s}(x)\\
 &= \int_{\mathcal{O}}\int_{\partial \mathcal{O}}\gamma(x)\frac{\partial }{\partial \mathbf{n}^a(x)}\mathcal{G}_\lambda(y,x) \d \mathrm{s}(x)\psi(y)\d y\\
 &= \int_{\mathcal{O}}\int_{\partial \mathcal{O}}\gamma(y)\frac{\partial }{\partial \mathbf{n}^a(y)}\mathcal{G}_\lambda(x,y) \d \mathrm{s}(y)\psi(x)\d x.
\end{align*}
To see \eqref{E57} note that 
\begin{align*}
\int_{\mathcal{O}}(\lambda- A) S(t)D_\lambda \gamma(x)\psi(x)\d x &= 
\int_{\mathcal{O}}D_\lambda \gamma(x)\left( \lambda - {{A}}^*\right)S^*(t)\psi(x)\d x \\
&=\int_{\partial \mathcal{O}}\gamma(x) \frac{\partial }{\partial \mathbf{n}^a}\left(A^*-\lambda \right)^{-1}\left( \lambda - {A}^*\right)S^*(t)\psi(x)\d \mathrm{s}(x)\\
&= - \int_{\partial \mathcal{O}}\gamma(x) \frac{\partial }{\partial {\mathbf{n}^a}} S^*(t)\psi(x)\d \mathrm{s}(x)\\
&=- \int_{\partial \mathcal{O}}\gamma(x) \frac{\partial }{\partial {\mathbf{n}^a}(x)} \int_{\mathcal{O}}G(t,y,x)\psi(y)\d y \d \mathrm{s}(x)\\
&= -\int_{\mathcal{O}}\int_{\partial \mathcal{O}} \gamma(x) \frac{\partial }{\partial {\mathbf{n}^a}(x)} G(t,y,x)\d \mathrm{s}(x) \psi(y)\d y \\
&=- \int_{\mathcal{O}}\int_{\partial \mathcal{O}} \gamma(y) \frac{\partial }{\partial {\mathbf{n}^a}(y)} G(t,x,y)\d \mathrm{s}(y) \psi(x)\d x. 
\end{align*}
\end{proof}
Note that by Assumption \ref{Ass61}, the coefficients $a_{i,j}$ are bounded. Therefore we have the following consequence of the theorem above.
\begin{corollary}\label{C64}
Under Assumptions \ref{A31}, \ref{A32}, and \ref{Ass61}, for any $T>0$ there is a constant $C>0$, such that for $t\in (0,T]$, $\psi\in L^2(\partial {\mathcal{O}},\d \mathrm{s})$, and $x\in \mathcal{O}$,  
$$
\left\vert S(t)B\psi(x)\right\vert = \left\vert (\lambda -A)S(t)D_\lambda \psi(x)\right\vert \le \frac{C}{\sqrt{t}}\left\vert \int_{\partial \mathcal{O}}g_{ct}(x-y) \psi(y)\d \mathrm{s}(y)\right\vert . 
$$
\end{corollary}
\section{Stochastic integration in $L^p$-spaces}\label{SRadonifying}

In this paper we need only very naive theory of stochastic integration in $L^p$-spaces. Namely, for $B:=(\lambda-A)D_\lambda$ set 
$$
\psi_k(t,x):= (S(t)Be_k)(x), \qquad x\in \mathcal{O}. 
$$
By Theorem \ref{T53}, 
$$
\psi_k(t,x)= -\int_{\partial \mathcal{O}}  \frac{\partial }{\partial {\mathbf{n}^a}(y)} G(t,x,y)e_k(y)\d \mathrm{s}(y). 
$$
We need to define 
$$
M(t,x):= \sum_{k} \int_0^t\psi_k(t-s,x)\d W_k(s),\quad x\in\mathcal{O},\,\,t\in[0,T]\,,
$$
and 
$$
M_\alpha (t,x):= \sum_{k} \int_0^t(t-s)^{-\alpha} \psi_k(t-s,x)\d W_k(s),\quad x\in\mathcal{O},\,\,t\in[0,T]\,,
$$
as $L^p_{\theta,\delta}$-valued provesses. 

The following  result from \cite{Brzezniak-Veraar} stated there as Proposition A.1, enables us to define rigorously each component of the sums above.  Below $W$ is a real valued Wiener process defined on a filtered probability space $(\Omega,\mathfrak{F}, (\mathfrak{F}_t), \mathbb{P})$. 
\begin{proposition}
Let $(\mathfrak{O}, \mathfrak{G}, \nu)$ be a $\sigma$-finite measurable space. Let $p,q\in (1,+\infty)$, $T\in (0,+\infty)$. For any adapted and strongly measurable process $\phi\colon [0,T]\times \Omega\mapsto L^p(\mathfrak{O})$ the following there assertions are equivalent.
\begin{enumerate}
\item There exists a sequence of adapted step processes $(\phi_n)$ such that 
\begin{align*}
&\lim_{n\to +\infty} \|\phi-\phi_n\|_{L^q(\Omega, L^p(\mathfrak{O}, L^2(0,T)))}=0, \\
&\left( \int_0^T\phi_n(t)\d W(t)\right)\quad \text{is a Cauchy sequence in $L^q(\Omega;L^p(\mathfrak{O}))$.}
\end{align*}
\item There exists a random variable $\eta \in L^q(\Omega;L^p(\mathfrak{O}))$ such that for all sets $A\in \mathfrak{G}$ with finite measure one has $(t,\omega)\mapsto \int_{A} \phi(t,\omega)\d \nu \in L^q(\Omega; L^2(0,T))$, and 
$$
\int_A \eta \d \nu = \int_0^T \int_{A} \phi(t)\d \nu \d W(t)\quad \text{in $L^q(\Omega)$.}
$$
\item $\|\phi\|_{L^q(\Omega;L^p(\mathfrak{O};L^2(0,T)))}<+\infty$. 
\end{enumerate}
Moreover, in this situation one has 
$$
\lim_{N\to +\infty} \int_0^T \phi_n(t)\d W(t)=\eta, 
$$
and there is a constant $C_{p,q}\in (0,+\infty)$ such that 
$$
C_{p,q}^{-1} \|\phi\|_{L^q(\Omega;L^p(\mathfrak{O};L^2(0,T)))} \le \|\eta \|_{L^q(\Omega;L^p(\mathfrak{O}))}\le C_{p,q}  \|\phi\|_{L^q(\Omega;L^p(\mathfrak{O};L^2(0,T)))}. 
$$
\end{proposition}

Process $\phi$ which satisfies ony of these conditions is called {\em $L^q$-stochastically integrable in $L^p(\mathfrak{O})$ on $[0,T]$} and we write 
$$
\int_0^T \phi(t)\d W(t):= \eta. 
$$

Given $\alpha \ge 0$ let 
\begin{equation} \label{E61}
\begin{aligned}
 \mathcal{J}_{T, \alpha}(\{e_k\},p, \theta, \delta) &:= \int_{\mathcal{O}}\left( \sum_{k} \int_0^T t^{-\alpha}\psi^2 _k(t,x)\d t\right)^{p/2} w_{\theta,\delta}(x)\d x. 
 \end{aligned}
 \end{equation}
By the  Burkholder--Davis--Gundy inequality, for every $p\in[1,+\infty)$ there exist positive constants $c_p$ and $C_p$ such that for all $t\in[0,T]$ 
$$
c_p \mathcal{J}_{T,0}(\{e_k\}, p,\theta, \delta)\le\mathbb{E}\int_{\mathcal{O}} \left\vert M(t,x)\right\vert ^p w_{\theta,\delta}(x)\d x \le C_p  \mathcal{J}_{T,0}(\{e_k\},p, \theta,\delta) \,,
$$
see for example \cite{veraar1}. We have thus the following result. 
 \begin{proposition}\label{P61}
 Given $T$, $p$, $\theta$ and $\delta$, the process
 $$
\sum_k \int_0^TS(T-t)Be_k\d W_k(t)
 $$
 takes values in $L^p_{\theta,\delta}$ if and only if $\mathcal{J}_{T,0}(\{e_k\}, p, \theta, \delta) <+\infty$. Moreover, if for a certain $\alpha>0$, $\mathcal{J}_{T,\alpha }(\{e_k\}, p, \theta, \delta) <+\infty$, then the process has continuous trajectories in $L^p_{\theta,\delta}$. 
 \end{proposition}

The simple idea above can be made rigorous and much more general (see e.g. \cite{Brzezniak-Neerven, Brzezniak-Veraar}). 
\section{Examples}\label{SExamples}
Let in the whole section Assumptions  \ref{A31}, \ref{A32}, and \ref{Ass61} be satisfied. 

\subsection{One dimensional case} Consider the simplest cases of $\mathcal{O}=(0,1)$ and $\mathcal{O}=(0,+\infty)$.  In the first case  the surface measure $\mathrm{s}= \delta_0+\delta_1$  and $L^2(\partial \mathcal{O}, \d \mathrm{s})\equiv \mathbb{R}^2$, whereas in the second case $\mathrm{s}= \delta_0$ and  $L^2(\partial \mathcal{O}, \d \mathbf{s})\equiv \mathbb{R}^1$. 
\begin{proposition}\label{P71}
Let  $p\in (1,+\infty)$ and $\theta \in (p-1,2p-1)$. Then the boundary problem 
$$
\begin{aligned}
\frac{\partial X}{\partial t}(t,x)&= \mathcal{A}X(t,x),\quad x\in (0,1),\\
X(t,0)&= \frac{\d W_0}{\d t}(t),\\
X(t,1)&= \frac{\d W_1}{\d t}(t),  
\end{aligned}
$$
defines Markov family  with continuous trajectories in the space $L^p_{\theta,0}$.  
\end{proposition}
\begin{proof} We are in a framework of Theorem \ref{T19}. By Theorem \ref{T34} the heat semigroup can be extended to $L^p_{\theta,0}$ for $\theta \in [0,2p-1)$. Clearly $H_W$ is 2-dimensional with $e_1=\chi_{\{0\}}$ and $e_2=\chi_{\{1\}}$.  Taking into account Proposition  \ref{P61}, it is enough to  verify whether for $\theta \in (p-1,2p-1)$ and $T\in (0, +\infty)$,   there is an $\alpha >0$ such that $\mathcal{J}:= \mathcal{J}_{T,\alpha }(\chi_{0},\chi_{1}, p, \theta, 0)<+\infty$. Let $\alpha >0$. By  Corollary \ref{C64}, we have 
\begin{align*}
\mathcal{J} &\le c_1\int_0^1 \left[ \int_0^T  t^{-1-\alpha } \left( g^2_{ct}(x)+ g^2_{ct}(x-1)\right)\d t  \right]^{p/2} \min\{ x^\theta, (1-x)^\theta\}\d x \\
&\le c_2 \int_0^1 \left[ \int_0^T t^{-1-\alpha } g^2_{ct}(t)(x)\d t \right]^{p/2} x^\theta \d x\\
&\le c_3 \int_0^1 \left[ \int_0^T t^{-2-\alpha } \e^{-\frac{x^2}{2ct}} \d t \right]^{p/2} x^\theta \d x= c_3 \int_0^1 \left[ \int_0^{T/x^2}  t^{-2-\alpha } \e^{-\frac{1}{2ct}}   \d t x^{-2-2\alpha} \right]^{p/2} x^\theta \d x\\
&\le c_3 \int_0^1 x^{\theta -p- \alpha p}\d x. 
\end{align*}
\end{proof}

Similar calculation can be done in the case of half-line $\mathcal{O}=(0,+\infty)$. 
\begin{proposition}\label{P72}
Assume that $\delta>1/2$ and  $p\in (1,+\infty)$ and $\theta \in (p-1,2p-1)$.   Let $\mathcal{A}$ be a second order defined as above. Then the boundary problem 
$$
\begin{aligned}
\frac{\partial X}{\partial t}(t,x)&= \mathcal{A}X(t,x),\quad x\in (0,+\infty),\\
X(t,0)&= \frac{\d W_0}{\d t}(t), 
\end{aligned}
$$
defines Markov family  with continuous trajectories in the space $L^p_{\theta,\delta }$. 
\end{proposition}
\begin{proof} We have 
\begin{align*}
\int_{0}^{+\infty} \left[ \int_0^T t^{-2-\alpha }\e^{-\frac{x^2}{ct}} \d t \right] ^{p/2} x^\theta (1+|x|^2)^{-\delta} \d x &\le c_1 \int_0^{+\infty} x^{\theta-p-\alpha p}  (1+|x|^2)^{-\delta} \d x.  
\end{align*}
\end{proof}
\begin{remark}
{\rm Assume that $A$ is equal to the Laplace operator $\Delta$.  Let $p\in (1,+\infty)$ and $\theta \in (p-1,2p-1)$. Then the Markov family defined by boundary problem on $(0,1)$ or $(0,+\infty)$ has a unique invariant measure. For, in the case of interval we can use Theorem \ref{T42}, whereas  in the case of problem on half-line  we can use a direct approach.  Namely, we have, see Remark \ref{Rem1} of Section \ref{SHalf-space}, 
$$
\frac{\partial G}{\partial \mathbf{n}_y}(t,x,0)= -\frac{x}tg_{2t}(x)= -\frac{x}{2t\sqrt{\pi t}}\e^{-\frac{x^2}{4t}}. 
$$
Clearly, there are  constant $C, c$ such that for all $t,x\ge 0$, 
$$
\left\vert \frac{\partial G}{\partial \mathbf{n}_y}(t,x,0)\right\vert \le \frac{C}{\sqrt{t}} g_{ct}(x). 
$$
Therefore, from the proof of the proposition it follows that $\sup_{T>0} \mathcal{J}_{T}(\chi_{0}, p,\theta,0)<+\infty$, and the desired conclusion holds, see Theorem \ref{T19}. }
\end{remark}

\subsection{Equation on a ball} Let $\mathcal{O}=\mathbb{B}^d$ be a unite ball in $\mathbb{R}^d$ with center at $0$. We assume that $d\ge 2$, otherwise we have  the case of equations on an interval studied in the previous section. Then $\partial \mathcal{O}= \mathbb{S}^{d-1}$.   Assume that the boundary noise has the form 
\begin{equation}\label{E71}
W(t,y)= \sum_{k} e_k(y)W_k(t), \qquad t\ge 0, \ y\in \mathbb{S}^{d-1}, 
\end{equation}
where $(e_k)$ is a sequence of  functions on $\mathbb{S}^{d-1}$ and $(W_k)$ are independent real-valued Wiener processes. 
\begin{proposition} \label{P74} If 
\begin{equation}\label{E72}
A:= \sum_{k} \sup_{y\in \mathbb{S}^{d-1}}e_k^2(y)<+\infty. 
\end{equation}
then for any $1<p$ and $\theta \in (p-1,2p-1)$, boundary problem \eqref{E11}  defines a Markov family with continuous trajectories in  $L^p_{\theta,0}$. 
\end{proposition}
\begin{proof} Let $\mathcal{J}:= \mathcal{J}_{T,\alpha }(\{e_k\},p,\theta,0)$. Using, Corollary \ref{C64}, \eqref{E61}, \eqref{E72},   we obtain  
\begin{equation}\label{E73}
\begin{aligned}
\mathcal{J} &\le c_1  \int_{\mathbb{B}^d } \left[\sum_{k} \int_0^T  t^{-1-\alpha} \left(  \int_{\mathbb{S}^{d-1}}  g_{ct}(x-y) \vert e_k(y)\vert    \d \mathrm{s}(y) \right)^2 \d t\right]^{p/2} w_{\theta,0}(x)\d x\\
&\le c_1A ^{p/2}  \int_{\mathbb{B}^d} \left[   \int_0^T  t^{-1-\alpha}  \left( \int_{\mathbb{S}^{d-1}} g_{ct}(x-y)   \d \mathrm{s}(y) \right)^2 \d t  \right]^{p/2} w_{\theta,0}(x)\d x\\
&\le c_2 \int_{\mathbb{B}^d} \left[   \int_0^T t^{-1-\alpha -d}   \left( \int_{\mathbb{S}^{d-1}} \e^{ -\frac{|x-y|^2}{ct} } \d \mathrm{s}(y) \right)^2 \d t  \right]^{p/2} w_{\theta,0}(x)\d x. 
\end{aligned}
\end{equation}
We have to evaluate 
$$
I(t,x):= \int_{\mathbb{S}^{d-1}}\e^{-\frac{|x-y|^2}{ct}} \d s(y). 
$$
Let $\rho(x):=\textrm{dist}\left(x,\mathbb{S}^{d-1}\right)$. We are showing that there is a  constant $C_1>0$ such that 
\begin{equation}\label{E74}
I(t,x)\le C_1 \e^{-\frac{\rho^2(x)}{C_1t}} t^{\frac{d-1}{2}}, \qquad \forall\, t\in (0,T], \ x\in \mathbb{B}^d. 
\end{equation}
In fact our method  leads also to the following lower bound estimate 
$$
C_2 \e^{-\frac{\rho^2}{C_3t}} t^{\frac{d-1}{2}}\le  I(t,x). 
$$
Our proof of \eqref{E74} is elementary. In Lemma \ref{L710}  from Section \ref{S73} we will establish estimate \eqref{E74} for an arbitrary bounded region in $\mathbb{R}^d$. In the proof of Lemma \ref{L710} we use the Laplace method. To see \eqref{E74}, note that we may assume that $x=(x_1,0,0,\ldots,0)= (1-\rho(x), 0,\ldots,0)$. Note that for $\rho,r\in [0,1]$ we have 
$$
1+\sqrt{1-r^2}-\rho \ge \left\vert 1-\sqrt{1-r^2}-\rho\right\vert. 
$$
Therefore 
\begin{align*}
I(t,x)&= \int_{\mathbb{S}^{d-1}}\exp\left\{ -\frac{\left\vert 1-\rho(x) - y_1\right\vert ^2 +\sum_{j=2}^d y_j^2}{ct}\right\} \d \mathrm{s}(y)
\\
&\le 2 \int_{\mathbb{B}^{d-1}} \exp\left\{ -\frac{\left\vert 1-\rho(x)  - \sqrt{1- \vert z\vert ^2}\right\vert  ^2 +\vert z\vert ^2 }{ct}\right\} \left( 1+ \sum_{j=1}^{d-1} \frac{z_j^2}{1-\vert z\vert ^2}\right)^{1/2}\d z\\
&\le c_3\int_0^1  \exp\left\{ -\frac{\left\vert 1-\rho(x)  - \sqrt{1- r ^2}\right\vert  ^2 +r ^2 }{ct}\right\} \left( 1-r^2\right)^{-1/2} r^{d-2} \d r.
\end{align*}

Note that  that for all $r\in [0,1]$ we have 
$$
\frac{1}{2} r^2 \le 1-\sqrt{1-r^2}\le  r^2.
$$
For, let $f(r)= (1-\sqrt{1-r})/r$. Then $f(0)=1/2$, $f(1)= 1$ and 
$$
f'(r)= \frac{(\sqrt{1-r}-1)^2}{2r^2 \sqrt{1-r}}\ge 0. 
$$
Therefore  there is a $c_4>0$ such that for all $r,\rho \in [0,1]$ we have 
$$
\left\vert 1-\rho  - \sqrt{1- r ^2}\right\vert  ^2 +r ^2 \ge c_4( \rho^2 + r^2). 
$$
For, we have
\begin{align*}
\left\vert 1-\rho  - \sqrt{1- r ^2}\right\vert  ^2 +r ^2&= \rho^2 + (1-\sqrt{1-r^2})^2 -2\rho(1-\sqrt{1-r^2})+r^2\\
&\ge \rho^2 + \frac{r^4}{4} - 2\rho r^2+r^2\\
&\ge \rho^2 + \frac{r^4}{4} - \rho^2 \kappa- \frac{r^4}{\kappa}+r^2= \rho^2(1-\kappa) + r^2\left(1- \frac{r^2}{\kappa}+ \frac{r^2}{4}\right). 
\end{align*}
Let $\kappa \in (4/5,1)$. Then 
$$
\left\vert 1-\rho  - \sqrt{1- r ^2}\right\vert  ^2 +r ^2\ge (1-\kappa)\rho^2 + r^2\left(\frac{5}{4} - \frac{1}{\kappa}\right),
$$
which gives the desired estimate. 

Summing up we have 

\begin{align*}
I(t,x)&\le c_3 \int_0^{1}  \exp\left\{ -\frac{c_4( \rho^2(x) + r^2)}{ct}\right\} \left(1-r^2\right)^{-1/2} r^{d-2} \d r \le \e^{-\frac{c_4\rho^2(x)}{ct}} C(t), 
\end{align*}
where 
\begin{align*}
C(t)&= c_3\int_0^{1}  \exp\left\{ -\frac{c_4 r^2}{ct}\right\} \left(1-r^2\right)^{-1/2} r^{d-2} \d r\\
&\le c_3\e^{-\frac{c_4}{4ct}}\int_{1/2}^1 \left(1-r^2\right)^{-1/2} \d r + c_3\frac{2}{\sqrt{5}} \int_0^{1/2} \exp\left\{ -\frac{c_4 r^2}{ct}\right\}  r^{d-2} \d r\\
&\le C_1 t^{\frac{d-1}{2}},
\end{align*}
which gives \eqref{E74}. 

Combining \eqref{E73} with \eqref{E74} we obtain 
\begin{align*}
\mathcal{J}&\le c_2  \int_{\mathbb{B}^d} \left[   \int_0^T t^{-1-\alpha -d}   \left( \int_{\mathbb{S}^{d-1}} \e^{ -\frac{|x-y|^2}{ct} } \d \mathrm{s}(y) \right)^2 \d t  \right]^{p/2} w_{\theta,0}(x)\d x\\
&\le c_5 \int_{\mathbb{B}^d} \left[   \int_0^T t^{-2-\alpha }    \e^{ -\frac{2 \rho^2(x) }{C_1t} } \d t  \right]^{p/2} w_{\theta,0}(x)\d x\\
&\le c_6 \int_{\mathbb{B}^d} \rho^{-p-\alpha p+\theta}(x) \d x\le  c_7 \int_0^1 r^{-p+\theta -\alpha p}\d r.  
\end{align*}
\end{proof}

\begin{remark} 
{\rm Assume that $A$ is equal to the Laplace operator $\Delta$.  Since the semigroup is exponentially stable on $L^p$, $1<p<+\infty $, then by  Theorem \ref{T42}, it is exponentially stable on $L^p_{\theta,\delta}$. Therefore the Markov family defined by the boundary problem \eqref{E11} on $L^p_{\theta,0}$ for $p-1<\theta <2p-1$, $p\in (1,+\infty)$ has a unique invariant measure.}
\end{remark}

Assumption \eqref{E72} ensures that $W$ is a random field on $[0,+\infty)\times \mathbb{S}^{d-1}$.  Below we present a natural example of the process satisfying the above assumptions. 
\begin{example}
Let $(W_{{k}})$ and $(\tilde W_{{k}})$  be sequences of independent Wiener  processes. Let $(a_{k})$ and $(b_k)$ be  sequences  of real numbers such that $\sum_{{k}}a_{{k}}^2<+\infty$. Then  
$$
W(t,y)= \sum_{{k}}  a_{{k}} \left(   \cos \langle y, b_{k}\rangle W_{k}(t)+ \sin \langle y, b_{k}\rangle\tilde  W_{k}(t)\right)
$$
 can obviously be written in the form \eqref{E71}. Moreover, condition \eqref{E72} is satisfied. 

For each $t\ge 0$, $W(t,\cdot) $ is rotational invariant random field on $\mathbb{S}^{d-1}$.  Indeed  we have  
\begin{align*}
\mathbb{E} \, W(t,y)W(t,z)& = \sum_{{k}} a_{k}^2 t \left( \cos  \langle y, b_{k}\rangle \cos  \langle z, b_{k}\rangle+ 
\sin  \langle y, b_{k}\rangle \sin   \langle z, b_{k}\rangle\right)\\
&=\sum_{{k}} a_{k}^2   \cos \langle b_k , y-z \rangle  t. 
\end{align*}
\end{example}
 
 In the case of the so-called \emph{white noise on $\mathbb{S}^{d-1}$}, $W$ is formally defined by $\eqref{E71}$ with $\{e_k\}$ being an orthonormal basis of $L^2(\mathbb{S}^{d-1}, \d \mathrm{s})$.
\begin{proposition}\label{P78}
Assume that $W$ is a white noise on $\mathbb{S}^{1}$. Let $p>1$ and $\theta\in  \left(\frac{3p}{2}-1,2p-1\right)$.  Then the boundary problem \eqref{E11}  defines a Markov family with continuous trajectories in  $L^p_{\theta,0}$.  If $A$ is equal to Laplace operator, then the Markov family defined by the boundary problem \eqref{E11} on $L^p_{\theta,0}$ for $p\in (1,+\infty)$ and $\theta \in \left(\frac{3p}{2}-1, 2p-1\right)$ has a unique invariant measure.
\end{proposition} 
\begin{proof} Assume that $d>1$. Let $\mathcal{J}:= \mathcal{J}_{T,\alpha}(\{e_k\},p,\theta,0)$. Using \eqref{E61}  and then Corollary \ref{C64},  we obtain  
\begin{align*}
\mathcal{J} &=  \int_{\mathbb{B}^d } \left[\sum_{k} \int_0^T t^{-\alpha} \left(  -\int_{\mathbb{S}^{d-1}}  \frac{\partial }{\partial {\mathbf{n}^a}(y)} G(t,x,y)e_k(y)\d \mathrm{s}(y) \right)^2 \d t\right]^{p/2} w_{\theta,0}(x)\d x\\
&= \int_{\mathbb{B}^d } \left[  \int_0^T t^{-\alpha}  \int_{\mathbb{S}^{d-1}}  \left\vert \frac{\partial }{\partial {\mathbf{n}^a}(y)} G(t,x,y)\right\vert ^2 \d \mathrm{s}(y)  \d t\right]^{p/2} w_{\theta,0}(x)\d x\\
&\le c_1 \int_{\mathbb{B}^d} \left[   \int_0^T   t^{-1-\alpha }\int_{\mathbb{S}^{d-1}} g^2_{ct}(x-y)   \d \mathrm{s}(y)  \d t  \right]^{p/2} w_{\theta,0}(x)\d x\\
&\le c_2 \int_{\mathbb{B}^d} \left[   \int_0^T t^{-1-d-\alpha }    \int_{\mathbb{S}^{d-1}} \e^{ -\frac{2|x-y|^2}{ct} } \d \mathrm{s}(y)  \d t  \right]^{p/2} w_{\theta,0}(x)\d x. 
\end{align*}
By  \eqref{E74} there is a constant $C>0$ such that  
\begin{align*}
\mathcal{J}&\le C\int_{\mathbb{B}^d} \left[   \int_0^T t^{-1-\alpha -d+ \frac{d-1}{2}}   \e^{ -\frac{\rho^2(x)}{Ct} }   \d t  \right]^{p/2} w_{\theta,0}(x)\d x\\
&\le C\int_{\mathbb{B}^d} \left[   \int_0^T t^{-\frac{d+3}{2}-\alpha }   \e^{ -\frac{\rho^2(x)}{Ct} }   \d t  \right]^{p/2} w_{\theta,0}(x)\d x\\
&\le C_1 \int_{\mathbb{B}^d} \rho^{\frac{p}{2}( 2-(d+3))-p\alpha}(x)w_{\theta,0}(x)\d x \le C_2 \int_0^1 r^{-\frac{p}{2}(d+1)+ \theta -p\alpha}\d r. 
\end{align*}
Therefore $\mathcal{J}<+\infty$ if $-1< -\frac{p}{2}(d+1)+ \theta -p\alpha$. In general case we need $\theta <2p-1$ and $p>1$. This requires 
$$
-1+\frac{p}{2}(d+1)+ p\alpha<\theta <2p-1\quad \text{and}\quad p>1. 
$$
Consequently, it is required that $d=1$ or $d=2$.  However the case $d=1$ has been already excluded.  
\end{proof}

\begin{remark} \label{R88}{\rm Let us drop the assumption that $e_k\in L^2(\partial \mathcal{O}, \d \mathrm{s}(y))$. Namely, assume that  $W(t,x)= B(t)\delta_{\widehat y}$, where $B$ is a real valued Wiener process and $\delta_{\widehat y}$ is the Dirac delta function at $\widehat y\in \mathbb{S}^{1}$. Then,  using \eqref{E61}  and then Corollary \ref{C64},  we obtain  
\begin{align*}
\mathcal{J}_{T,\alpha}(\delta_{\widehat y} ,p,\theta,0)&=  \int_{\mathbb{B}^2 } \left[\sum_{k} \int_0^T t^{-\alpha} \left(   \frac{\partial }{\partial {\mathbf{n}^a}(y)} G(t,x,\widehat {y}) \right)^2 \d t\right]^{p/2} w_{\theta,0}(x)\d x\\
&\le c_1 \int_{\mathbb{B}^2} \left[   \int_0^T t^{-1-2-\alpha}    \e^{ -\frac{2\vert x-\widehat y\vert ^2}{ct} } \ \d t  \right]^{p/2} w_{\theta,0}(x)\d x\\
&\le c_2 \int_{\mathbb{B}^2}   \vert x-\widehat y\vert ^{-2p-\alpha p}\rho^\theta (x)\d x\le c_3 \int_0^1 r^{-2p+ \theta -\alpha p}\d r. 
\end{align*}
We need $-2p+\theta -\alpha p>-1$, which is in contradiction with  $\theta <2p-1$. Therefore, we cannot treat this case. }
\end{remark}
\subsection{The case of a bounded region in $\mathbb{R}^d$}\label{S73} Let $\mathcal{O}$ be a bounded $C^{1,\alpha}$ region in $\mathbb{R}^d$, $d\ge 2$.  Set  
$$
I(t,x):= \int_{\partial \mathcal{O}}\e^{-\frac{|x-y|^2}{ct}} \d s(y), \qquad t\in (0,T],  \ x\in \mathcal{O}. 
$$
Recall that $\rho(x):=\textrm{dist}\left(x,\partial \mathcal{O}\right)$. We have the following generalization of \eqref{E74} established for $\mathcal{O}=\mathbb{B}^d$. 
\begin{lemma}\label{L710}
There exist  constants $C_1,C_2>0$ such that 
$$
I(t,x)\le C_1 t^{\frac{d-1}{2}}\e^{-\frac{\rho^2(x)}{C_2t}}, \qquad \forall\, t\in (0,T], \ x\in \mathcal{O}. 
$$
\end{lemma}
\begin{proof}
Fix $\ve>0$ and $x\in\mathcal{O}$. Let $\varepsilon>0$ be fixed. Since $\mathcal{O}$ is a bounded $C^{1,\alpha}$-domain, there exist open sets $\mathcal{O}_i\subset\R^d$, $i=1,2,\ldots,n$, such that for every $i$:
\begin{enumerate}
\item
for every $i$ there exists, up to a shift and a rotation, a $C^{1,\alpha}$ function $a_i$ such that $x\in\mathcal{O}\cap\mathcal{O}_i$ if and only if $x=\la \bar x,x_d-a_i\la\bar x\ra\ra$ with $x_d-a_i\la\bar x\ra>0$. 
\item
$\partial\mathcal{O}\subset\bigcup_i\mathcal{O}_i$\,,
\item 
if $x\in\mathcal{O}\setminus \bigcup_i\mathcal{O}_i$ then $\rho(x)>\varepsilon$,
\item
for each $i\le n$ there exists a $C^{1,\alpha}$-diffeomorphism 
\[g^i\colon \mathcal{O}_i\to C^d:=\left\{z\in\R^d\colon \,-1<z_k<1,\,k=1,\ldots d\right\}\] such that $h^i:=\la g^i\ra^{-1}\colon C^d\to\mathcal{O}_i$ is of class $C^{1,\alpha}$ as well and 
\[g^i\la\mathcal{O}_i\cap\partial\mathcal{O}\ra=C^{d-1}=\left\{z\in C^d\colon z_d=0\right\},\]
\[\mathcal{O}_i\cap\mathcal{O}=\left\{z\in C^d\colon z_d>0\right\},\]
\[\mathcal{O}_i\cap\overline{\mathcal{O}}^c=\left\{z\in C^d\colon z_d<0\right\}.\]
\end{enumerate}
Assume that $x\in\mathcal{O}_i\cap\mathcal{O}$ for certain $i\le n$. Then for $y=\la \bar y,0\ra\in\mathcal{O}_i\cap\partial\mathcal{O}$ 
\[|x-y|^2=\left|\la\bar x,x_d-a_i\la\bar x\ra\ra-\la\bar y,0\ra\right|^2\ge\left|\bar x-\bar y\right|^2+\rho^2(x)\,,\]
hence
\[\begin{aligned}
I(t,x)&\le Ct^{\frac{d-1}{2}}\e^{-\frac{\rho^2(x)}{ct}}\int_{\mathcal{O}_i\cap\partial\mathcal{O}}g_{ct}\la\bar x-\bar y\ra\d s(y)\\
&= Ct^{\frac{d-1}{2}}\e^{-\frac{\rho^2(x)}{ct}}\int_{C^{d-1}}g_{ct}\la h^i(u)-h^i(v)\ra J\la h^i\ra(v)\d v\,.
\end{aligned}\]
Since $h^i\colon C^d\to\mathcal{O}_i$, there exists a constant $c_1>0$ such that 
\[\left|h^i(u)-h^i(v)\right|\ge c_1|u-v|,\quad u,v\in C^d\,.\]
 Therefore,  
\[\begin{aligned}
I(t,x)&\le Ct^{\frac{d-1}{2}}\e^{-\frac{\rho^2(x)}{ct}}\int_{C^{d-1}}g_{ct}(u-v)\d v \le Ct^{\frac{d-1}{2}}\e^{-\frac{\rho^2(x)}{ct}}\,,
\end{aligned}\]
and the lemma follows for $x\in\mathcal{O}\cap\mathcal{O}_i$ for every $i\le n$.

If $x\in\mathcal{O}\setminus\bigcup_i\mathcal{O}_i$, then $\rho(x)\ge\ve$ and the lemma trivially follows. 
\end{proof}
Recall that 
$$
W(t,y)= \sum_{k} e_k(y)W_k(t), \qquad t\ge 0, \ y\in \partial \mathcal{O}, 
$$
where $(e_k)$ is a sequence of  functions on $\partial \mathcal{O}$ and $(W_k)$ are independent real-valued Wiener processes. 
\begin{proposition}\label{P711}
$(i)$ Assume that $\sum_{k} \sup_{y\in \partial \mathcal {O}}e_k^2(y)<+\infty$.  Then for any $1<p$ and $\theta \in (p-1,2p-1)$, the boundary problem \eqref{E11}  defines a Markov family with continuous trajectories in   $L^p_{\theta,0}$. 

\noindent
$(ii)$ Assume that  $d=2$ and $W$ is  a white-noise on $\partial \mathcal{O}$. Then for $1<p$ and $\theta \in \left(\frac{3p}{2}-1,2p-1\right)$ the boundary problem \eqref{E11}  defines a Markov family with continuous trajectories in  $L^p_{\theta,0}$. 
\end{proposition}
\begin{proof}[Proof of $(i)$] Using the calculations from the proof of Proposition \ref{P74}, and then our Lemma \ref{L710} we obtain 
\begin{align*}
 \mathcal{J}_{T,\alpha}(\{e_k\},p,\theta,0) &\le c_1 \int_{\mathcal{O}} \left[   \int_0^T t^{-1-d-\alpha }   \left( I(t,x) \right)^2 \d t  \right]^{p/2} w_{\theta,0}(x)\d x\\
 &\le c_2  \int_{\mathcal{O}} \left[   \int_0^T t^{-1-d+d-1-\alpha }    \e^{-2\frac{\rho^2(x)}{C_1t}}  \d t  \right]^{p/2} w_{\theta,0}(x)\d x\\
 &\le c_3 \int_{\mathcal{O}} \rho(x) ^{-p+\theta-p\alpha}(x)\d x. 
\end{align*}
It is easy to show, see the proof of Lemma \ref{L45}, that the integral is finite if and only if $-p+\theta -p\alpha >-1$. 

\end{proof}
\begin{proof}[Proof of $(ii)$] Using the calculations from the proof of Proposition \ref{P78}, and then our Lemma \ref{L710} we obtain 
\begin{align*}
 \mathcal{J}_{T,\alpha}(\{e_k\},p,\theta,0) &\le c_1 \int_{\mathcal{O}} \left[   \int_0^T t^{-1-d-\alpha}    \int_{\partial \mathcal{O}} \e^{ -\frac{2|x-y|^2}{ct} } \d \mathrm{s}(y)  \d t  \right]^{p/2} w_{\theta,0}(x)\d x \\
 &\le c_2  \int_{\mathcal{O}} \left[   \int_0^T t^{-1-d+\frac{d-1}{2}-\alpha }    \e^{-2\frac{\rho^2(x)}{C_1t}}  \d t  \right]^{p/2} w_{\theta,0}(x)\d x\\
 &\le c_3 \int_{\mathcal{O}} \rho^{-\frac{p}{2} (d+1) +\theta-\alpha p}(x)\d x. 
\end{align*}

\end{proof}

\subsection{Half-space with spatially homogeneous Wiener process}\label{SHalf-space}
In this section $\mathcal{O}=(0,+\infty)\times \mathbb{R}^m$, and $W$ is the so-called spatially homogeneous Wiener process on $\mathbb{R}^m\equiv \{0\}\times \mathbb{R}^m= \partial \mathcal{O}$. We adopt the notation $x=(x_0,x_1,\ldots, x_n)=(x_0,\mathbf{x})$. 

\begin{definition} 
A process $W$ taking values in the space of tempered distributions $\mathcal{S}'(\mathbb{R}^m)$ is called a \emph{spatially homogeneous Wiener process} if and only if: 
\begin{itemize}
\item[$(i)$] It is  Gaussian process with continuous trajectories in $\mathcal{S}^\prime(\mathbb{R}^m)$.
\item[$(ii)$] For each  $\psi\in\mathcal{S}(\mathbb{R}^m)$, $t\mapsto  \left( {W}(t),\psi\right)$ is a one  dimensional Wiener process.  
\item[$(iii)$] For each fixed $t\ge 0$ the law of $W(t)$  is invariant with  respect to all translations
$\tau ^\prime_h\colon \mathcal{S}^\prime(\mathbb{R}^m)\rightarrow \mathcal{S}^\prime(\mathbb{R}^m)$,
$h\in \mathbb{R}^m$, where  $\tau _h\colon \mathcal{S}(\mathbb{R}^m)\rightarrow \mathcal{S}(\mathbb{R}^m)$,
$\tau _h \psi (\cdot) = \psi (\cdot +h)$ for $\psi \in \mathcal{S}(\mathbb{R}^m)$.  
\end{itemize}
\end{definition}

The law of  a spatially homogeneous Wiener process  $W$ on $\mathbb{R}^m$ is characterized by its \emph{spectral measure} $\mu$ on $\left(\mathbb{R}^m,\mathcal{B}(\mathbb{R}^m)\right)$. Recall, see \cite{Peszat-Zabczyk1} that $\mu$ is a positive symmetric Radon tempered measure on $\mathbb{R}^m$, and for any test functions $\psi,\phi\in \mathcal{S}(\mathbb{R}^m)$, 
$$
\mathbb{E}\langle W(t),\psi\rangle \langle W(s),\phi\rangle = t\wedge s\int_{\mathbb{R}^m} \mathcal{F}\psi(\mathbf{x})\overline{\mathcal{F}\phi(\mathbf{x})}\mu(\d \mathbf{x}), 
$$
where $\mathcal{F}$ denotes the Fourier transform.  Note, see e.g. \cite{Peszat-Zabczyk1} that if the spectral measure is finite, then $W$ is a random field, such that for any $\mathbf{x}$, $W(\cdot,\mathbf{x})$ is a one dimensional Wiener process. Moreover, for fixed $t$, the field $W(t,\mathbf{x})$ is stationary in $\mathbf{x}$. 

Let 
$$
L^2_{(s)}(\mu):=\left\{ \psi\in  L^2(\mathbb{R}^d\mapsto\mathbb{C},\mathcal{B}(\mathbb{R}^m),\mu)\colon \psi(-\mathbf{x})=\overline{\psi(\mathbf{x})}\right\}. 
$$
Then, see  \cite{Peszat-Zabczyk1},  the Reproducing Kernel Hilbert Space $H_W$ of $W$ is given by 
$$
H_W=\left\{\mathcal{F}\psi\colon \psi\in L^2_{(s)}(\mu)\right\}
$$
and 
$$
\langle \mathcal{F}\psi, \mathcal{F}\phi\rangle_{H_W}= \int_{\mathbb{R}^d}\psi(\mathbf{x})\overline {\phi(\mathbf{x})}\mu(\d \mathbf{x}). 
$$
Thus any  orthonormal basis $\{e_k\}$ of $H_W$ is of the form  $e_k= \mathcal{F}(f_k\mu)$, where $\{f_k\}$ is an orthonormal basis of $L^2_{(s)}(\mu)$.   We will identify $\mathbb{R}^m$ with $\partial (0,+\infty)\times \mathbb{R}^m=\{0\}\times \mathbb{R}^m$.  

We have
\begin{align}\nonumber
\mathcal{J}_{T,\alpha}(\{e_k\},p,\theta,\delta)&:= \int_{\mathcal{O}} \left[ t^{-\alpha}\int_0^T \sum_{k} \left( \int_{\mathbb{R}^m}  \frac{\partial G}{\partial \mathbf{n}_{\mathbf{y}} }(t,x,(0,\mathbf{y})) e_k(\mathbf{y})\d \mathbf{y}  \right)^2 \d t \right]^{p/2}  w_{\theta,\delta}(x)\d x \\ \nonumber
&= \int_{\mathcal{O}} \left[ \int_0^T t^{-\alpha} \sum_{k} \left( \int_{\mathbb{R}^m}  \frac{\partial G}{\partial \mathbf{n}_{\mathbf{y}} }(t,x,(0,\mathbf{y})) \mathcal{F}(f_k\mu)(\mathbf{y})\d \mathbf{y} \right)^2 \d t \right]^{p/2}  w_{\theta,\delta}(x)\d x \\ \label{E77}
&\le \int_{\mathcal{O}} \left[ \int_0^T t^{-\alpha} \int_{\mathbb{R}^m}  \left\vert \mathcal{F}^{-1}_{\mathbf{y}} \frac{\partial G}{\partial \mathbf{n}_{\mathbf{y}} }(t,x,(0,\mathbf{y}))\right\vert ^2 \mu(\d \mathbf{y}) \d t \right]^{p/2}  w_{\theta,\delta}(x)\d x. 
\end{align}
\begin{proposition}\label{P713}
Assume that:  the spectral measure of $W$ is finite,   $\delta>(m+1)/2$, $p\in (1,+\infty)$ and $\theta \in (p-1,2p-1)$. Then   boundary problem \eqref{E11} defines Markov family  with continuous trajectories in  $L^p_{\theta,\delta }$. 
\end{proposition}
\begin{proof}
If the measure $\mu$ is finite, then
\begin{align*}
 \int_{\mathbb{R}^m}  \left\vert \mathcal{F}^{-1}_{\mathbf{y}} \frac{\partial G}{\partial \mathbf{n}_{\mathbf{y}} }(t,x,(0,\mathbf{y}))\right\vert ^2 \mu(\d \mathbf{y})
 &\le \mu(\mathbb{R}^m) \sup_{\mathbf{y}\in \mathbb{R}^m} \left\vert \mathcal{F}^{-1}_{\mathbf{y}} \frac{\partial G}{\partial \mathbf{n}_{\mathbf{y}} }(t,x,(0,\mathbf{y}))\right\vert ^2\\
 &\le \mu(\mathbb{R}^m) \left[ \int_{\mathbb{R}^m}  \left\vert  \frac{\partial G}{\partial \mathbf{n}_{\mathbf{y}} }(t,x,(0,\mathbf{y}))\right\vert \d \mathbf{y}\right]^2 \\
 &\le c t^{-2} \e ^{-\frac {x_0^2}{ct}}, 
\end{align*}
where in the last estimate we use \eqref{E61}.  We have 
$$
\int_0^{T} t^{-2-\alpha} \e ^{-\frac {x_0^2}{ct}}\d t \le \int_0^{+\infty} s^{-2-\alpha } \e^{-\frac{1}{2s}}\d s x_0^{-2-2\alpha}\le c_1 x_0^{-2}.
$$
Therefore, by \eqref{E77},   we have 
\begin{align*}
\mathcal{J}_{T,\alpha}(\{e_k\},p,\theta,\delta)&\le c_1 \int_0^{+\infty}\int_{\mathbb{R}^m}   \min\{x_0,1\}^{\theta} x_0^{-p-\alpha p} \left( 1+ \vert x_0\vert ^2+\vert \mathbf{x}\vert ^2\right)^{-\delta} \d x_0 \d \mathbf{x} \\
&\le c_2+ c_2 \int_0^1 x_0^{\theta-p-\alpha p}\d x_0. 
\end{align*}
\end{proof}

\begin{remark}\label{Rem1}
{\rm If $\mathcal {A}=\Delta$ then $G$ is given by \eqref{EK1}. Then, with $\overline {x}$ defined by \eqref{EK2}, 
$$
\frac{\partial G}{\partial \mathbf{n}_y}(t,x,y)= -\frac{\partial G}{\partial y_0}(t,x,y)= \frac{y_0-x_0}{2t}g_{2t}(x-y)+ \frac{-y_0-x_0}{2t}g_{2t}(\overline x-y).
$$
}
\end{remark}

\begin{remark}
{\rm  \eqref{E61} gives estimates for   $\frac{\partial G}{\partial \mathbf{n}_y}(t,x,(0,\mathbf{y}))$. Unfortunately, we are not able to use them to compare the Fourier transforms of  $\frac{\partial G}{\partial \mathbf{n}_y}(t,x,(0,\mathbf{y}))$ and $\tilde g_{ct}$. This problem can be solved under a certain technical assumption on the spectral measure $\mu$, see the lemma below. For a measure for which this assumption is violated see \cite{Peszat}. }
\end{remark}

Write 
$$
K_\alpha(r):=\int_0^{+\infty} s^{-2-\alpha}\e^{-\frac{1}{s}-r^2 s}\d s,\qquad r\ge 0, 
$$

\begin{lemma} \label{L716}Assume that either $A=\Delta$ or there is a finite symmetric measure $\mu_0$ on $\mathbb{R}^m$ such that  $\mathcal{F}\mu+ \mathcal{F}\mu_0$ is a non-negative  measure. Let $\delta> (m+1)/2$, $p>1$ and $\theta\in(p-1,2p-1)$.  Then for any $T>0$, there is a constant $C>0$ such that 
$$
\mathcal{J}_{T,\alpha}(\{e_k\},p,\theta,\delta)\le C+C\int_{\mathbb{R}^m}\int_0^{+\infty} \left[ \int_{\mathbb{R}^m} K_\alpha\left( \frac{x_0\vert \mathbf{y}\vert}{C}\right)  \mu(\d \mathbf{y})\right]^{p/2} x_0^{-p-\alpha p}w_{\theta,\delta}((x_0,\mathbf{y}))\d x_0\d \mathbf{y}.
$$
\end{lemma}
\begin{proof} Taking into account \eqref{E77} and Proposition \ref{P713} we can assume that $\mathcal{F}\mu$ is a non-negative measure.  Then 
\begin{align*}
& \int_{\mathbb{R}^m}  \left\vert \mathcal{F}^{-1}_{\mathbf{y}} \frac{\partial G}{\partial \mathbf{n}_{\mathbf{y}} }(t,x,(0,\mathbf{y}))\right\vert ^2 \mu(\d \mathbf{y})\\
 &= \int_{\mathbb{R}^m} \int_{\mathbb{R}^m} \frac{\partial G}{\partial \mathbf{n}_{\mathbf{u}} }(t,x,(0,\mathbf{u}))\frac{\partial G}{\partial \mathbf{n}_{\mathbf{v}} }(t,x,(0,\mathbf{v}))\mathcal{F}\mu({\mathbf{u}} -{\mathbf{v}} )\d {\mathbf{u}} \d {\mathbf{v}} \\
&\le c_1 ^2 t^{-1} \int_{\mathbb{R}^m} \int_{\mathbb{R}^m} g_{ct}(x-(0,\mathbf{u}))g_{ct}(x- (0,\mathbf{v}))\mathcal{F}\mu({\mathbf{u}} -{\mathbf{v}} )\d {\mathbf{u}} \d {\mathbf{v}} \\
&\le c_2t^{-2} \e^{-\frac{x_0^2}{c_2t}}\int_{\mathbb{R}^m}  \e^{-\frac{t\vert \mathbf{z}\vert ^2}{c_2}} \mu(\d \mathbf{z}). 
\end{align*}
Therefore, by \eqref{E77}, 
\begin{align*}
\mathcal{J}_{T,\alpha}(\{e_k\},p,\theta,\delta)&\le \int_{\mathbb{R}^m}\int_0^{+\infty} \left[ \int_{\mathbb{R}^m} \int_0^T c_2t^{-2-\alpha} \e^{-\frac{x_0^2}{c_2t}}  \e^{-\frac{t\vert \mathbf{z}\vert ^2}{c_2}} \d t \mu(\d \mathbf{z})\right]^{p/2} w_{\theta,\delta}((x_0,\mathbf{y}))\d x_0\d \mathbf{y}.
\end{align*}
Since 
$$
\int_0^T c_2t^{-2-\alpha} \e^{-\frac{x_0^2}{c_2t}}  \e^{-\frac{t\vert \mathbf{z}\vert ^2}{c_2}} \d t \le x_0^{-2-2\alpha} \int_0^{+\infty} s^{-2} \e^{-\frac{1}{s} -\frac{x_0^2\vert\mathbf{z}\vert ^2 s }{c_2^2}} \d s= x_0^{-2-\alpha} K_\alpha\left(\frac{x_0\vert \mathbf{z}\vert}{c_2}\right), 
$$
the desire conclusion follows. 
\end{proof}
If $\mu(\d \mathbf{y})=\d \mathbf{y}$, then $W$ is the co-called \emph{cylindrical Wiener process on} $L^2(\mathbb{R}^m)$ or equivalently $\frac{\partial W}{\partial t}(t,\mathbf{y})$, $t\ge 0$, $\mathbf{y}\in \mathbb{R}^m$,   is the \emph{space-time white noise}. Then $\mathcal{F}\mu$ is the Dirac delta measure. We have the following consequence of Lemma \ref{L716} and our general Theorem \ref{T19}.  Note that our  results on heat semigroups on weighted  spaces do not allow $m>1$. 
\begin{proposition} \label{P717}Let $m=1$, $\delta >1$, and let $W$ be  a cylindrical Wiener process on $L^2(\mathbb{R})$. Then the boundary problem \eqref{E11} defines a Markov family  with continuous  trajectories in the space $L^p_{\theta,\delta }$ for $p>1$ and $\theta\in \left(\frac{3p}{2}-1,2p-1\right)$. 
\end{proposition}
\begin{proof}
We have 
\begin{align*}
\int_{\mathbb{R}^m}  K_\alpha \left(\frac{x_0\vert \mathbf{z}\vert}{C}\right)  \mu(\d \mathbf{z})&= \int_0^{+\infty} s^{-2-\alpha }\e^{-\frac{1}{s}}\int_{\mathbb{R}^m}\e^{-\frac{x_0^2  \vert \mathbf{z}\vert ^2}{C^2}s }\d \mathbf{z}\d s\\
&= \int_0^{+\infty} s^{-2-\alpha}\e^{-\frac{1}{s}}\left( 2\pi \frac{C^2}{2x_0^2s } \right)^{m/2}\d s= \tilde Cx_0^{-m}. 
\end{align*}
Then, by Lemma \ref{L716},   $\mathcal{J}_{T,\alpha}(\{e_k\},p,\theta,\delta)$ is finite if 
$$
\int_0^1 r^{-\frac{mp}{2}-p +\theta-\alpha p}\d r<+\infty. 
$$
This requires $\theta >-1+ \frac{mp}{2} +p$.  Since we need $2p-1>\theta$ we arrive at the condition $p>\frac{mp}{2}$.  Since $p\ge 1$ we are able to deal with the boundary problem only in the case of $m=1$.
\end{proof}

In many interesting cases, see e.g. \cite{Dalang}, the spectral measure $\mu$ of $W$ is absolutely continuous and its density is the so-called  \emph{Bessel potential}. Namely for a parameter $\kappa>0$ let 
$$
\mu_\kappa(\d \mathbf{y}) :=  \left(1+\vert \mathbf{y}\vert ^2\right)^{-\kappa/2}\d \mathbf{y}.
$$
Then the space correlation $\Gamma_\kappa (\mathbf{y}):= \mathcal{F} \mu_\kappa(\mathbf{y})$  is a non-negative continuous function on $\mathbb{R}^m\setminus\{0\}$. Moreover,  asymptotically as $\vert \mathbf{y}\vert \to 0$,
\begin{equation}\label{E78}
\Gamma_\kappa (\mathbf{y}) \approx  \begin{cases} C_{m,\kappa}\vert \mathbf{y} \vert ^{\kappa-m}&\text{for $0<\kappa<m$,}\\
C_{m,\kappa}\log \frac{1}{\vert \mathbf{y}\vert} &\text{for $\kappa=m$,}\\
C_{m,\kappa}&\text{for $\kappa>m$.}
\end{cases}
\end{equation}
Finally 
\begin{equation}\label{E79}
\Gamma_\kappa (\mathbf{y})  \approx C_{m,\kappa}\e ^{-\vert \mathbf{y}\vert }, \qquad \text{as $\vert \mathbf{y}\vert \to +\infty$}.
\end{equation}

Note that in the limit $\kappa\downarrow 0$ we obtain Lebesgue measure corresponding to the cylindrical Wiener process on $L^2(\mathbb{R}^m)$ treated in Proposition \ref{P717}. 
\begin{proposition}\label{P718}
Let $W$ be a  Wiener process on $\mathbb{R}^m$  with the spectral measure $\mu_\kappa$, $0<\kappa\le m$.  $(i)$ If $\kappa \ge m$, then boundary problem \eqref{E11} defines a Markov family with continuous trajectories   in the space $L^p_{\theta,\delta }$ for $p>1$ and $\theta \in (p-1,2p-1)$.

\noindent  $(ii)$ If $m-2<\kappa <m$, then boundary problem \eqref{E11} defines a Markov family  with continuous trajectories in the space $L^p_{\theta,\delta }$ for 
$$
1< p< +\infty, \qquad p+\frac{p}{2}(m-\kappa)-1<\theta <2p-1. 
$$
\end{proposition}
\begin{proof} Note that if $\kappa >m$ then $\mu_\kappa$ is finite and we may apply our Proposition \ref{P713}.  Therefore we restrict our attention to the case of $0<\kappa \le m$.  We have 
\begin{align*}
\int_{\mathbb{R}^m}  K_\alpha\left(\frac{x_0\vert \mathbf{z}\vert}{C}\right)  \mu_\kappa(\d \mathbf{z})&= \int_0^{+\infty} s^{-2-\alpha }\e^{-\frac{1}{s}}\int_{\mathbb{R}^m}\e^{-\frac{x_0^2 \vert \mathbf{z}\vert ^2}{C}s}\left(1+\vert \mathbf{z}\vert ^2\right)^{-\kappa/2}  \d \mathbf{z}\, \d s\\
&= \int_0^{+\infty} s^{-2-\alpha }\e^{-\frac{1}{s}}\int_{\mathbb{R}^m}\left[ \mathcal{F}_{\mathbf{z}}^{-1} \e^{-\frac{x_0^2 s }{C}\vert \mathbf{z}\vert ^2} \right] \Gamma_\kappa(\mathbf{z}) \d \mathbf{z}\, \d s\\
&= \int_0^{+\infty} s^{-2-\alpha }\e^{-\frac{1}{s}}\int_{\mathbb{R}^m}  \left( 2\pi \sigma^2 \right)^{-m/2} \e^{-\frac{\vert \mathbf{z}\vert ^2}{2\sigma^2}}\Gamma_\kappa(\mathbf{z}) \d \mathbf{z}\,\d s,
\end{align*}
where   $\sigma^2 :=  \frac{x_0^2s }{2C}$. Then by \eqref{E79}, 
$$
\int_{\{ \vert \mathbf{z}\vert \ge 1\}} (2\pi \sigma^2)^{-m/2} \e^{-\frac{\vert \mathbf{z}\vert ^2}{2\sigma^2}}\Gamma_\kappa(\mathbf{z}) \d \mathbf{z}\le  C\int_{\{ \vert \mathbf{z}\vert \ge 1\}} (2\pi \sigma^2)^{-m/2} \e^{-\frac{\vert \mathbf{z}\vert ^2}{2\sigma^2}- \frac{\vert \mathbf{z}\vert }{C}} \d \mathbf{z}\le C_1<+\infty,
$$
where $C_1<+\infty$ does not depend on $\sigma$. 

Next, by \eqref{E78}, if $\kappa =m$, then 
\begin{align*}
\int_{\{ \vert \mathbf{z}\vert < 1\}} (2\pi \sigma^2)^{-m/2} \e^{-\frac{\vert \mathbf{z}\vert ^2}{2\sigma^2}}\Gamma_\kappa(\mathbf{z}) \d \mathbf{z}&\le  C\int_{\{ \vert \mathbf{z}\vert < 1\}} (2\pi \sigma^2)^{-m/2} \e^{-\frac{\vert \mathbf{z}\vert ^2}{2\sigma^2}} \log \frac{1}{C \vert  \mathbf{z}\vert} \d \mathbf{z}\\
&\le C_1 \int_0^1 (2\pi \sigma^2)^{-m/2} \e^{-\frac{r ^2}{2\sigma^2}} \left\vert \log C r\right\vert  r^{m-1} \d r.
\end{align*}
If $m>1$, then $r\mapsto \left\vert \log C r\right\vert  r^{m-1}$ is a continuous function on a closed interval $[0,1]$. Therefore 
$$
\int_{\{ \vert \mathbf{z}\vert < 1\}} (2\pi \sigma^2)^{-m/2} \e^{-\frac{\vert \mathbf{z}\vert ^2}{2\sigma^2}}\Gamma_1(\mathbf{z}) \d \mathbf{z}\le \tilde C_1, 
$$
where again $\tilde C_1<+\infty$ does not depend on $\sigma$. 

If $\kappa=m=1$ then $r\mapsto \vert \log Cr\vert$ in integrable on $[0,1]$ with any power $>1$. Therefore, using the H\"older inequality we obtain that for any $q>1$ there is an independent of $\sigma$  constant $C(q)$ such that 
$$
\int_{\{ \vert \mathbf{z}\vert < 1\}} (2\pi \sigma^2)^{-1/2} \e^{-\frac{\vert \mathbf{z}\vert ^2}{2\sigma^2}}\Gamma_1(\mathbf{z}) \d \mathbf{z}\le 
C(q)\left( \int_0^1 (2\pi \sigma^2)^{-q/2} \e^{-q\frac{r ^2}{2\sigma^2}}  \d r\right)^{1/q}\le C(q)\sigma^{\frac{1-q}{2q}}. 
$$
Taking  $q>1$  small enough we conclude that for any $\varepsilon>0$ there is an independent of $\sigma$  constant $C_\varepsilon$ such that 
$$
\int_{\{ \vert \mathbf{z}\vert < 1\}} (2\pi \sigma^2)^{-1/2} \e^{-\frac{\vert \mathbf{z}\vert ^2}{2\sigma^2}}\Gamma_1(\mathbf{z}) \d \mathbf{z}\le 
 C_\varepsilon \sigma^{-\varepsilon}. 
$$

Finally if $0<\kappa <m$, then 
\begin{align*}
\int_{\{ \vert \mathbf{z}\vert < 1\}} (2\pi \sigma^2)^{-m/2} \e^{-\frac{\vert \mathbf{z}\vert ^2}{2\sigma^2}}\Gamma_\kappa(\mathbf{z}) \d \mathbf{z}&\le 
C_2\int_0^1 (2\pi \sigma^2)^{-m/2} \e^{-\frac{r ^2}{2\sigma^2}}  r^{\kappa-m+m-1}  \d r\\
&\le  C_2\int_0^1 (2\pi \sigma^2)^{-m/2} \e^{-\frac{r ^2}{2\sigma^2}}  r^{\kappa-1}  \d r\\
&\le C_3 \sigma ^{-m+1+\kappa -1} \int_0^{+\infty}   \e^{-\frac{u ^2}{2}}  u^{\kappa-1}  \d r\\
&\le C_4\sigma^{\kappa -m}. 
\end{align*}
Summing up, we see that there is an independent of $x_0$  constant $c_1$ such that 
$$
\int_{\mathbb{R}^m}  K_\alpha \left(\frac{x_0\vert \mathbf{z}\vert}{C}\right)  \mu_\kappa(\d \mathbf{z}) \le c_1,\qquad \text{if $\kappa=m>1$,}
$$
and 
$$
\int_{\mathbb{R}^m}  K_\alpha \left(\frac{x_0\vert \mathbf{z}\vert}{C}\right)  \mu_\kappa(\d \mathbf{z}) \le c_1+ c_1 x_0^{\kappa-m},\qquad \text{if $\kappa<m$.}
$$
If $\kappa=m=1$ then for any $\varepsilon >0$ then  there is an independent of $x_0$  constant $c(\varepsilon)$ such that 
$$
\int_{\mathbb{R}^m}  K_\alpha\left(\frac{x_0\vert \mathbf{z}\vert}{C}\right)  \mu_\kappa(\d \mathbf{z}) \le c_1 + c(\varepsilon)x_0^{-\varepsilon}.
$$

Thus, by Lemma \ref{L716},   $\mathcal{J}_{T,\alpha}(\{e_k\},p,\theta,\delta)$ is finite if 
$$
\int_{0}^1 x_0^{-p+\theta-p\alpha} d x_0 <+\infty,\qquad \text{if $\kappa=m\ge 1$}
$$
and 
$$
\int_{0}^1x_0 ^{-p+ \frac{p}{2}(\kappa -m)+\theta-p\alpha } \d x_0<+\infty,\qquad \text{if $\kappa<m$.}
$$
\end{proof}
\begin{remark}
{\rm Note that,  we are not able to treat the case of $m> 2$ and $\kappa \le m-2$.}
\end{remark}
\appendix 
\section{Proof of Lemma \ref{L45}}\label{AppL}
Assume that $\mathcal{O}$ is a half space. Without any loss of generality we may assume that $\mathcal{O}=\{x=(x_1,\mathbf{x})\in \mathbb{R}^d\colon x_1>0\}$. We also assume that $d>1$.  Then 
\begin{align*}
\int_{\mathcal{O}}g_{ct}(x-y)\rho^\alpha(y)\, \d y&= \int_0^{+\infty} \int_{\mathbb{R}^{d-1}} g_{ct}(x-y)y^\alpha\,   \d y\\
&= \int_0^{+\infty} \left(2\pi ct\right)^{-1/2} \e^{-\frac{(x_1-y_1)^2}{2ct}} y^\alpha _1\, \d y_1\\
&\le t^{\frac{\alpha}{2}}  \int_{\mathbb{R}} \left(2\pi c\right)^{-1/2} \e^{-\frac{(t^{-\frac 12}x_1-z)^2}{2c}} \vert z\vert ^{\alpha} \d z.
\end{align*}
Since 
\begin{align*}
&\sup_{r\in \mathbb{R} } \int_{\mathbb{R}} (2\pi c)^{-1/2}\e^{-\frac{(r-z)^2}{2c}} \vert z\vert ^{\alpha} \d z\\
&\le \int_{\vert z\vert \le 1} (2\pi c)^{-1/2}\vert z\vert ^{\alpha} \d z+\sup_{r\in \mathbb{R}} \int_{\mathbb{R}} (2\pi c)^{-1/2}\e^{-\frac{(z-r)^2}{2c}} \d z  <+\infty,
\end{align*}
the desired conclusion follows. \qquad $\square$

\bigskip 
Assume that  $\mathcal{O}$ is a bounded $C^{1,\alpha}$-domain. For $y\in\mathbb{R}^d$ we will write $y=\left(y_1,\mathbf{y}\right)\in\mathbb{R}^{d-1}\times\mathbb{R}$ and $g_{ct}(y)=g_{ct}^{(1)}\left( y_1\right) g_{ct}^{(d-1)}\left(\mathbf{y}\right)$. Since $\mathcal  O$ is a bounded $C^{1,\alpha}$-domain, its boundary $\partial\mathcal{O}$ can be covered with a finite number of open sets $\mathcal{O}_i$, such that for every $i$ there exists a $C^{1,\alpha}$ function $h^i$ such that (up to a shift and rotation of the domain) 
$$
\mathcal{O}\cap\mathcal{O}_i=\left\{y\in\mathbb{R}^d\colon  y_1>h^i\left(\mathbf{ y}\right)\right\}\,.
$$
Moreover, for $t$ small enough we have 
$$
\mathcal{O}_t\subset\bigcup_i\mathcal{O}_i\,.
$$
For each $i$ we can define a $C^1$-diffeomorphism 
$$
g^i\colon \mathcal{O}_i\to g^i\left(\mathcal{O}_i\right),\quad g^i\left(y_1, \mathbf{y}\right)=\left(y_1-h^i\left(\mathbf{y}\right), \mathbf{y}\right)
$$
Clearly, the Jacobian $J^i$ of $g^i$ satisfies the condition $\left\vert J^i(x)\right\vert =1$. Therefore, for $z=\left( z_1, \mathbf{z}\right)\in g^i\left(\mathcal{O}_i\right)$ we have 
$$
z_1=y_1-h^i\left(\mathbf{y}\right)=\inf_{v\in\mathcal { O}_i\cap\partial\mathcal{O}}\left\vert g^i\left( y\right)-g^i\left( v\right)\right\vert \,.
$$
Since 
$$
c_1\vert y-v\vert \le \left\vert g^i\left( y\right)-g^i\left( v\right)\right\vert \le c_2\vert y-v\vert ,\quad y,v\in\mathcal{O}_i\cap\mathcal{O}
$$
we find that 
$$
c_1 z_d\le \rho(y)=\rho\left(\bar y,y_d\right)\le c_2z_d,\quad y\in\mathcal{O}_t\,.
$$
We have 
\begin{align*}
\int_{\mathcal{O}_t}\rho^{\alpha}(y) g_{ct}(x-y)\d y&=\sum_{i}\int_{\mathcal  {O}_t\cap\mathcal {O}_i}\rho^{\alpha}(y) g_{ct}(x-y)\d y =:\sum_i J_i(t,x)\,,
\end{align*}
and it is enough to show that for every $i$ 
$$
\sup_{t\le 1}\sup_{x\in\mathcal{O}} t^{-\frac{\alpha}{2}}J_i(t,x) <+\infty\,.
$$
Changing variables  we obtain 
\begin{align*}
J_i(t,x)&\le C\int_{\mathbb{R}^{d}}\vert z_1\vert ^{\alpha} g_{ct} (z-x)\d z = Ct^{\frac \alpha 2} \int_{\mathbb{R}}\vert y \vert  ^{\alpha} g_{c}^{(1)}\left(x_dt^{-1/2} - y\right)\d y\\
&\le Ct^{\frac \alpha 2}\left(2\pi c\right) ^{-\frac 12}  \int_{-1}^1 \vert y \vert ^{\alpha}\d y + Ct^{\frac{\alpha} 2} \int_{\{ \vert y\vert \ge 1\}} g_{c}^{(1)}\left(x_dt^{-1/2} - y\right)\d y\\
&\le Ct^{\frac \alpha 2}\left(2\pi c\right) ^{-\frac 12}  \int_{-1}^1 \vert y \vert ^{\alpha}\d y + Ct^{\frac{\alpha} 2}. \qquad \square
\end{align*}

\section{$C_0$-property without Assumption \ref{Axx}}\label{AppT} We are showing that for any $p\in [1,+\infty)$ and $\theta\in[0,p)$ there exists  a constant $M_{p,\theta}$ such that 
\begin{equation}\label{E36}
\vert  S(t)\psi \vert _{L^p_\theta}\le M_{p,\theta} \vert \psi\vert _{L^p_\theta}, \qquad \forall\, t\in (0,1], \ \forall\, \psi \in L^p. 
\end{equation}
By \eqref{E33} and the Jensen inequality we have 
\begin{align*}
\left\vert S(t)\psi\right\vert ^p_{L^p_\theta}&\le C^p\int_{\mathcal  O}\rho^\theta (x)\int_{\mathcal  O}m_t^p(y)g_{ct}(x-y)|\psi(y)|^p\,\d y\,\d x.
\end{align*}
Changing  variables we obtain 
\begin{equation} \label{E37}
\begin{aligned}
&\left\vert S(t)\psi\right\vert ^p_{L^p_\theta}\\
&\qquad \le C^pt^{d/2}\int_{\mathcal  O/\sqrt{t}}\rho^\theta (x\sqrt{t})\int_{\mathcal  O/\sqrt{t}}m_t^p(y\sqrt{t})g_{c}(x-y)|\psi(y\sqrt{t})|^p\,\d y\,\d x\\
&\qquad \le C^pt^{d/2}\int_{\mathcal  O/\sqrt{t}}\rho^\theta (x\sqrt{t})\int_{\mathcal  O/\sqrt{t}}m_t^p(y\sqrt{t})g_{c}(x-y)|\psi(y\sqrt{t})|^p\,\d y\,\d x.
\end{aligned}
\end{equation}
Recall that $\rho(x)= \textrm{dist}\left(x,\partial \mathcal {O}\right)$.  Define 
 $$
 (\mathcal{O}/\sqrt{t})_1:= \left\{x\in \mathcal{O}/\sqrt{t} \colon m_t(x\sqrt{t})=\rho(x\sqrt{t})/\sqrt{t}<1\right\}
 $$
 and 
 $$
 (\mathcal{O}/\sqrt{t})_1^c:= \left\{x\in \mathcal{O}/\sqrt{t}\colon m_t(x\sqrt{t})= 1\right\}. 
 $$
 Then we have 
 $$
 \left\vert S(t)\psi\right\vert ^p_{L^p_\theta}\le C^pt^{d/2}\left( I_1+I_2+I _3+I_4\right), 
 $$
 where 
 \begin{align*}
 I_1&:= \int_{(\mathcal{O}/\sqrt{t})_1} \d x \, \rho^{\theta} (x\sqrt{t})\int_{(\mathcal{O}/\sqrt{t})_1}  \d y\, m_t^p (y\sqrt{t}) g_{c}(x-y)  |\psi(y\sqrt{t})|^{p},\\
 I_2&:=\int_{(\mathcal{O}/\sqrt{t})_1} \d x \, \rho^{\theta} (x\sqrt{t})\int_{(\mathcal{O}/\sqrt{t})_1^c} \d y  \, m_t^p(y\sqrt{t})g_{c}(x-y) |\psi(y\sqrt{t})|^{p},\\
 I_3&:=  \int_{(\mathcal{O}/\sqrt{t})_1^c} \d x\, \rho^{\theta} (x\sqrt{t})\int_{(\mathcal{O}/\sqrt{t})_1^c}  \d y \,  m_t^p(y\sqrt{t})g_{c}(x-y) |\psi(y\sqrt{t})|^{p},\\
 I_4&:= \int_{(\mathcal{O}/\sqrt{t})_1^c} \d x\, \rho^\theta  (x\sqrt{t})\int_{(\mathcal{O}/\sqrt{t})_1}  \d y\, m_t^p(y\sqrt{t}) g_{c}(x-y)  |\psi(y\sqrt{t})|^{p}. 
 \end{align*}
Set $\phi(y)=\psi(y\sqrt{t})$.  Taking into account that $\rho(z\sqrt{t})/ \sqrt{t}\le 1$ for $z\in (\mathcal{O}/\sqrt{t})_1$, and the fact that since $\theta <p$
$$
\frac{\rho^p (y\sqrt{t})}{t^{p/2}}\le \frac{\rho^\theta  (y\sqrt{t})}{t^{\theta /2}},\quad y\in (\mathcal{O}/\sqrt{t})_1,
$$
we have 
\begin{align*}
I_1 &=  \int_{(\mathcal{O}/\sqrt{t})_1} \d x \, \rho^{\theta} (x\sqrt{t})\int_{(\mathcal{O}/\sqrt{t})_1}  \d y\, \frac{\rho^p (y\sqrt{t})}{t^{p/2}} g_{c}(x-y)  |\phi(y)|^p\\
&\le  \int_{(\mathcal{O}/\sqrt{t})_1} \d x \, \frac{\rho^{\theta}(x\sqrt{t})}{t^{\theta /2}} \int_{(\mathcal{O}/\sqrt{t})_1}  \d y\,  \rho^{\theta } (y\sqrt{t})g_{c}(x-y)  |\phi(y)|^p\\
&\le  \int_{(\mathcal{O}/\sqrt{t})_1} g_{c}(x-y) \d x \int_{(\mathcal{O}/\sqrt{t})_1}  \d y\, \rho^{\theta} (y\sqrt{t})  |\phi(y)|^p\\
&\le  \int_{\mathbb{R}^d} g_{c}(x-y) \d x \int_{(\mathcal{O}/\sqrt{t})_1}  \d y\, \rho^{\theta} (y\sqrt{t})  |\phi(y)|^p\\
&\le C_1 \int_{(\mathcal{O}/\sqrt{t})_1}  \d y\, \rho^{\theta} (y\sqrt{t})  |\phi(y)|^p\le  C_1 \vert \phi\vert ^p_{L^p_\theta(\mathcal{O}/\sqrt{t})}. 
\end{align*}
Taking into account that $(\mathcal{O}/\sqrt{t})_1^c$,   $1\le \rho(z\sqrt{t})/ \sqrt{t}$ we find that 
 \begin{align*}
 I_2 &=   \int_{(\mathcal{ O}/\sqrt{t})_1} \d x \, \rho^{\theta} (x\sqrt{t}) \int_{(\mathcal{ O}/\sqrt{t})_1^c} \d y  \, g_{c}(x-y) |\phi(y)|^p \\
&\le \int_{\mathbb{R}^d} g_{c}(x-y)\d x t^{\theta /2}\int_{(\mathcal{ O}/\sqrt{t})_1^c} \d y   |\phi(y)|^p \\
&\le C_1\int_{(\mathcal{O}/\sqrt{t})_1^c} \d y  \frac{\rho^\theta (y\sqrt{t})}{t^{\theta /2}}  t^{\theta /2}|\phi(y)|^p \\
&\le C_1 \vert\phi\vert ^p_{L^p_\theta(\mathcal{O}/\sqrt{t})}
 \end{align*}
and 
 \begin{align*}
 I_3&=    \int_{(\mathcal{O}/\sqrt{t})_1^c} \d x\, \rho^{\theta} (x\sqrt{t}) \int_{(\mathcal{O}/\sqrt{t})_1^c}  \d y \,  g_{c}(x-y) |\phi(y)|^{p}\\
 &=   \int_{(\mathcal{O}/\sqrt{t})_1^c} \d x\, \frac{\rho^{\theta} (x\sqrt{t})}{\rho^{\theta} (y\sqrt{t})} g_{c}(x-y) \int_{(\mathcal{O}/\sqrt{t})_1^c}  \d y \, \rho^{\theta} (y\sqrt{t})  |\phi(y)|^{p}\\
 &\le  \sup_{y\in (\mathcal{O}/\sqrt{t})_1^c} \int_{(\mathcal{O}/\sqrt{t})_1^c} \d x\, \frac{\rho^{\theta} (x\sqrt{t})}{\rho^{\theta} (y\sqrt{t})} g_{c}(x-y) \vert \phi\vert ^p_{L^p_\theta(\mathcal{O}/\sqrt{t})}, 
\end{align*}
and finally, as  
$$
\frac{\rho^p(y\sqrt{t})}{t^{p/2}}\le \frac{\rho^\theta (y\sqrt{t})}{t^{\theta /2}}, \qquad y\in {(\mathcal{O}/\sqrt{t})_1}, 
$$
we have 
\begin{align*}
I_4&=   \int_{(\mathcal{O}/\sqrt{t})_1^c} \d x\, \rho^\theta  (x\sqrt{t}) \int_{(\mathcal{O}/\sqrt{t})_1}  \d y\,  \frac{\rho^p(y\sqrt{t})}{t^{p/2}} g_{c}(x-y)  |\phi(y)|^{p}\\
&\le  \int_{(\mathcal{O}/\sqrt{t})_1^c} \d x\, \rho^\theta  (x\sqrt{t}) \int_{(\mathcal{O}/\sqrt{t})_1}  \d y\,  \frac{\rho^\theta (y\sqrt{t})}{t^{\theta /2}} g_{c}(x-y)  |\phi(y)|^{p}\\
&\le  \int_{(\mathcal{O}/\sqrt{t})_1^c} \d x\, \frac{\rho^\theta  (y\sqrt{t})}{t^{\theta/2}} g_{c}(x-y)\int_{(\mathcal{O}/\sqrt{t})_1}  \d y\,  \rho^\theta (x\sqrt{t})   |\phi(y)|^{p}\\
&\le  \sup_{y \in (\mathcal{O}/\sqrt{t})_1}\int_{(\mathcal{O}/\sqrt{t})_1^c} \d x\, \frac{\rho^\theta  (y\sqrt{t})}{t^{\theta/2}} g_{c}(x-y)\vert \phi\vert ^p_{L^p_\theta(\mathcal{ O}/\sqrt{t})}. 
 \end{align*}
Note that 
$$
\vert \phi\vert ^p_{L^p_\theta(\mathcal{ O}/\sqrt{t})}= t^{-d/2}\vert \psi\vert ^p_{L^p_\theta}. 
$$
Therefore the proof will be completed as soon as we show that 
$$
A_1:= \sup_{t\in (0,1]}\sup_{y\in (\mathcal{O}/\sqrt{t})_1^c} \int_{(\mathcal{O}/\sqrt{t})_1^c}  \frac{\rho^{\theta} (x\sqrt{t})}{\rho^{\theta} (y\sqrt{t})} g_{c}(x-y)\d x <+\infty
$$
and 
$$
A_2:= \sup_{t\in (0,1]}\sup_{y \in (\mathcal{O}/\sqrt{t})_1}\int_{(\mathcal{O}/\sqrt{t})_1^c} \frac{\rho^\theta  (x\sqrt{t})}{t^{\theta/2}} g_{c}(x-y)\d x<+\infty. 
$$
To do this note that  
$$
\rho(z\sqrt{t})= \textrm{dist}\left(z\sqrt{t},\partial {\mathcal{O}}\right)= \sqrt{t}\, \textrm{dist}\left(z,\partial {\mathcal{O}/\sqrt{t}}\right). 
$$ 
Therefore 
$$
A_1= \sup_{t\in (0,1]}\sup_{y\in (\mathcal{O}/\sqrt{t})_1^c} \int_{(\mathcal{O}/\sqrt{t})_1^c}  \left(\frac{\textrm{dist}\left(x,\partial \mathcal{O}/\sqrt{t}\right)}{\textrm{dist}\left(y,\partial \mathcal{O}/\sqrt{t}\right)} \right)^\theta g_{c}(x-y)\d x
$$
and 
$$
A_2= \sup_{t\in (0,1]}\sup_{y\in (\mathcal{O}/\sqrt{t})_1} \int_{(\mathcal{O}/\sqrt{t})_1^c}  \textrm{dist}\left(x,\partial \mathcal{O}/\sqrt{t}\right)^\theta  g_{c}(x-y)\d x. 
$$
Given a domain $\mathcal{D}\subset \mathbb{R}^d$, $\mathcal{D}\not =\mathbb{R}^d$,  set 
\begin{align*}
\mathcal{D}_1&=\left\{x\in \mathcal{D}\colon \textrm{dist}\left(x,\partial \mathcal{D}\right)<1\right\},\\
\mathcal{D}_1^c&=\left\{x\in \mathcal{D}\colon \textrm{dist}\left(x,\partial \mathcal{D}\right)\ge 1\right\},\\
A_1(\mathcal{D})&= \sup_{y\in \mathcal{D}_1^c} \int_{\mathcal{D}_1^c}  \left(\frac{\textrm{dist}\left(x,\partial \mathcal{D}\right)}{\textrm{dist}\left(y,\partial \mathcal{D}\right)} \right)^\theta \exp\left\{ -\frac{|x-y|^2}{c}\right\}\d x, \\
A_2(\mathcal{D})&= \sup_{y\in \mathcal{D}_1} \int_{\mathcal{D}_1^c}  \textrm{dist}\left(x,\partial \mathcal{D}\right)^\theta \exp\left\{ -\frac{|x-y|^2}{c}\right\}\d x. 
\end{align*}
We have to show that there is a constant $N$ (independent of $\mathcal{D}$ but it can depend on $d$, $\theta$ and $c$) such that 
\begin{equation}\label{E38}
A_1(\mathcal{D})+ A_2(\mathcal{D})\le N. 
\end{equation}
We note first that for any $x,y\in\R^d$ 
$$
\left\vert \textrm{dist}\left(x,\partial \mathcal{D}\right)-\textrm{dist}\left(y,\partial \mathcal{D}\right)\right\vert \le |x-y|\,.
$$
We will consider $A_1(\mathcal D)$ first. For every $y\in\mathcal D_1^c$ we obtain 
\[\begin{aligned}
\left(\frac{\textrm{dist}\left(x,\partial D\right)}{\textrm{dist}\left(y,\partial\mathcal D\right)}\right)^\theta=&\left(\frac{\textrm{dist}\left(x,\partial\mathcal D\right)-\textrm{dist}\left(y,\partial\mathcal D\right)}{\textrm{dist}\left(y,\partial\mathcal D\right)}+1\right)^\theta\\
\le&\left(|x-y|+1\right)^\theta\,.
\end{aligned}\]
Therefore, 
$$
\begin{aligned}
A_1(\mathcal{D})\le&\sup_{y\in \mathcal{D}_1^c} \int_{\mathcal{D}_1^c}\left(|x-y|+1\right)^\theta\exp\left\{ -\frac{|x-y|^2}{c}\right\}\d x\\
\le&\sup_{y\in \mathcal{D}_1^c}\int_{\mathbb{R}^d}\left(|x-y|+1\right)^\theta\exp\left\{ -\frac{|x-y|^2}{c}\right\}\d x\\
=&\int_{\mathbb{R}^d}(1+|z|)^\theta\exp\left\{ -\frac{|z|^2}{c}\right\}\d x<+\infty. 
\end{aligned}
$$
Consider now $A_2(\mathcal D)$. Then, by similar arguments for every $y\in\mathcal D_1$ we obtain 
$$
\begin{aligned}
\textrm{dist}\left (x,\partial \mathcal{D}\right)^\theta\le&\left( \textrm{dist}\left(x,\partial \mathcal{D}\right)-\textrm{dist}\left(y,\partial \mathcal{D}\right)+\textrm{dist}\left(y,\partial \mathcal{D}\right)\right)^\theta
\\
\le&\la |x-y|+1\ra^\theta
\end{aligned}
$$
and again 
$$
A_2(\mathcal D)\le \int_{\mathbb{R}^d}(1+|z|)^\theta\exp\left\{ -\frac{|z|^2}{c}\right\}\d z<+\infty\,.
$$
Combining the two estimates above we obtain \eqref{E38} with 
$$
N=2\int_{\mathbb{R}^d}(1+|z|)^\theta\exp\left\{ -\frac{|z|^2}{c}\right\}\d z<+\infty\,.\qquad \square
$$

We are showing now the the gradient estimates.  Using Assumption  \ref{A32} (\eqref{E32} and \eqref{E33})  and the Jensen inequality we obtain 
$$
\begin{aligned}
\left| \frac{\partial S(t)\psi }{\partial x_j}\right|^p_{L^p_\theta}&=\int_{\mathcal{O}}\rho^\theta(x)\left|\int_{\mathcal{O}}\frac{\partial }{\partial x_j} G(t,x,y)\psi(y)dy\right|^p\d x\\
&\le \frac{C^p}{t^{p/2}}\int_{\mathcal{O}}\rho^\theta(x)\int_{\mathcal{O}}m_t^p(y)g_{ct}(x-y)|\psi(y)|^p\d y \d x\\
&=\frac{C^pt^{d/2}}{t^{p/2}}\int_{\mathcal{O}/\sqrt{t}}\rho^\theta(x\sqrt{t})\int_{\mathcal{O}/\sqrt{t}}m_t^p(y\sqrt{t})g_{c}(x-y)|\psi(y\sqrt{t})|^p\d y \d x\\
&=\frac{C^pt^{d/2}}{t^{p/2}}\left(I_1+I_2+I_3+I_4\right)\,,
\end{aligned}
$$
where $I_i$, $i=1, 2,3,4$ are defined in the previous section. Therefore we can use the estimates for $I_i$ and the desired conclusion follows. 
 $\qed$

\section*{Acknowledgment}
We would like to  thank Professor  Zdzislaw  Brze\'zniak  for very useful discussions on the topic.


\begin{thebibliography}{99}
\bibitem{Alos-Bonacorsi1}
E. Al\`os and  S. Bonaccorsi,
 {\em Stability for stochastic partial differential equations with Dirichlet white-noise boundary conditions},   Infin. Dimens. Anal. Quantum Probab. Relat. Top.  {\bf  5} (2002),  465--481.
 
\bibitem{Alos-Bonacorsi2}
 E. Al\`os and S. Bonaccorsi, {\em Stochastic partial differential equations with Dirichlet white-noise boundary conditions},  Ann. Inst. H. Poincar\'e Probab. Statist. {\bf 38} (2002), 125--154.

\bibitem{Balakrishnan}
A.V. Balakrishnan, {\em Applied Functional Analysis}, Springer-Verlag, Berlin Heidenberg New York, 1981.


\bibitem{bgpr}
Z. Brze\'zniak,  B. Goldys, S. Peszat,  and F. Russo, {\em Second order PDEs with Dirichlet white noise boundary conditions},  J. Evol. Equ. {\bf 15} (2015), 1--26.


\bibitem{Brzezniak-Neerven}
Z. Brze\'zniak and J. van Neerven, {\em Space-time regularity for linear stochastic evolution equations driven by spatially homogeneous noise}, J. Math. Kyoto Univ. 
{\bf 43}  (2003), 261--303. 


\bibitem{Brzezniak-Peszat}
Z. Brze\'zniak and S. Peszat, {\em  Hyperbolic equations with random boundary conditions}, in  Recent Development in Stochastic Dynamics and Stochastic Analysis,  (J. Duan, S. Luo and C. Wang, eds.)  World Scientific, 2010, pp. 1--22.

\bibitem{Brzezniak-Veraar}
Z.  Brze\'zniak and  M. Veraar, {\em Is the stochastic parabolicity condition dependent on $p$ and $q$?}, Electron. J. Probab.  {\bf 17}  (2012), 1--24.

\bibitem{Chi-Kim-Park}
S. Cho, P. Kim,  and H. Park, {\em Two-sided estimates on Dirichlet heat kernels for time-dependent parabolic operators with singular drifts in $C^{1,\alpha}$-domains},  J. Differential Equations {\bf 252} (2012), 1101--1145.


\bibitem{Dalang-Leveque} R. Dalang and O. L\'ev\^eque, {\em Second order linear hyperbolic SPDE's driven by isotropic Gaussian noise on a sphere}, Ann. Probab. {\bf 32}  (2004),  1068--1099.

\bibitem{Dalang}
R. C. Dalang, M. Santz-Sol\'e, \emph{H\"older--Sobolev regularity of the solution to the stochastic wave equation in dimension $3$},  Memoirs of the American Mathematical Society {\bf 199} (2006). 

\bibitem{DaPrato-Zabczyk}
G. Da Prato and J. Zabczyk, {\em Evolution equations with white-noise boundary conditions},  Stochastics Stochastics Rep. {\bf 42} (1993), 167--182.

\bibitem{DaPrato-Zabczyk2}
G. Da Prato and  J. Zabczyk, {\em  Stochastic Equations in Infinite Dimensions}, Second Edition,  Cambridge Univ. Press, Cambridge, 2014.

\bibitem{DaPrato-Zabczyk3}
G. Da Prato and  J. Zabczyk, {\em  Ergodicity for Infinite Dimensional Systems}, Cambridge Univ. Press, Cambridge, 1996.

\bibitem{Delfour}
A. Bensoussan, G. Da Prato, M.C. Delfour, and S.K. Mitter, {\em Representation and Control of Infinite Dimensional Systems (Systems \& Control: Foundations \& Applications)}, Birkh\''auser, Boston, 2006.

\bibitem{duncan}
 T.E. Duncan, B. Maslowski and B. Pasik-Duncan, B. {\em Ergodic boundary/point control of stochastic semilinear systems}. SIAM J. Control Optim. 36 (1998), no. 3, 1020--1047
\bibitem{Eidelman-Ivasishen}
S.D. Eidel'man and S.D. Ivasishen, {\em Investigations of the Green matrix for a homogeneous parabolic boundary value problem}, Trans. Moscow Math. Soc. {\bf 23} (1970), 179--242.

\bibitem{Fabri-Goldys}
G. Fabbri G. and B. Goldys, {\em  An LQ problem for the heat equation on the halfline with Dirichlet boundary control and noise} SIAM J. Control Optim. {\bf 48} (2009), 1473--1488.

\bibitem{Freidlin-Sowers}
M. Freidlin  and R. Sowers, {\em Central limit results for a reaction-diffusion equation with fast-oscillating boundary perturbations}, Stochastic partial differential equations and their applications (Charlotte, NC, 1991), pp. 101--112, Lecture Notes in Control and Inform. Sci., 176, Springer, Berlin, 1992.

\bibitem{Garling}
D.J.H. Garling, {\em Inequalities: A Journey into Linear Analysis},   Cambridge Univ.  Press,  2007. 

\bibitem{Krylov}
N.V. Krylov, {\em  The heat equation in $L^q((0,T),L^p)$-spaces with weights},  SIAM J. Math. Anal. {\bf 32} (2001), 1117--1141.

\bibitem{krylov_c1}
Kyeong-Hun Kim and N.V. Krylov, 
{\em On the Sobolev space theory of parabolic and elliptic equations in C1 domains}, SIAM J. Math. Anal. {\bf 36} (2004), 618--642. 

\bibitem{Lasiecka-Triggiani}
I. Lasiecka and R. Triggiani, {\em Control Theory for Partial Differential Equations Volume 1: Abstract Parabolic Systems}, Cambridge Univ. Press, Cambridge, 2000.

\bibitem{Lindemulder-Veraar}
N. Lindemulder and M. Veraar, {\em  The heat equation with rough boundary conditions and holomorphic functional calculus}, J.  Differential Equations {\bf 269} (2020) 5832--5899. 

\bibitem{Lions-Magenes}
J.L. Lions and E. Magenes, {\em Non-Homogeneous Boundary Value Problems and Applications I}, Springer-Verlag, Berlin Heidenberg New York, 1972.

\bibitem{Lunardi}
A. Lunardi, {\em Analytic Semigroups and Optimal Regularity in Parabolic Problems}, Birkhauser, 1995.

\bibitem{Maslowski}
B. Maslowski, {\em  Stability of semilinear equations with boundary and pointwise noise},   Ann. Scuola Norm. Sup. Pisa Cl. Sci. {\bf 22},  (1995), 55--93.

\bibitem{mazya}
 V. Maz'ya, \em{Sobolev Spaces with Applications to Elliptic Partial Differential Equations},  Second, revised and augmented edition,  Springer, 2011.

\bibitem{mclean}
W. McLean, \em{Strongly Elliptic Systems and Boundary Integral Equations}, Cambridge University Press, 2000.

\bibitem{munteanu}
 I. Munteanu, {\em Boundary stabilization of parabolic equations},  Birkh\"auser, Springer, Cham, 2019

\bibitem{munteanu_p}
I. Munteanu, {\em Stabilization of stochastic parabolic equations with boundary-noise and boundary-control}. J. Math. Anal. Appl. 449 (2017), 829--842
\bibitem{Mora}
X. Mora, {\em Semilinear parabolic problems define semiflows on $C^k$ spaces}, Trans. Amer. Math. Soc. {\bf 278}, (1983), 21--55.

\bibitem{Peszat}
S. Peszat, {\em The Cauchy problem for a nonlinear stochastic wave equation in any dimension}, J. Evol. Equ. {\bf 2} (2002), 383--394.

\bibitem{Peszat-Zabczyk1}
S. Peszat and J. Zabczyk, {\em Stochastic evolution equations with a spatially homogeneous Wiener process}, Stochastic Processes Appl. {\bf 72} (1997), 187--204.

\bibitem{Peszat-Zabczyk2}
S. Peszat and J. Zabczyk, {\em Stochastic Partial Differential Equations Driven by L\'evy Processes}, Cambridge Univ. Press, Cambridge, 2007.

\bibitem{Schwartz}
L.  Schwartz, {\em Th\'eorie des Distributions I, II},  Hermann \& Cie., Paris, 1950, 1951.

\bibitem{Solonnikov}
V.A. Solonnikov,  {\em  Green matrices for parabolic boundary value problems}, Sem. Math. V.A. Stieklov Math. Inst. Leningrad {\bf 14} (1969), 132--150.


\bibitem{Sowers}
R.B. Sowers, {\em  Multidimensional reaction-diffusion equations with white noise boundary perturbations},   Ann. Probab.  {\bf 22} (1994), 2071--2121.

\bibitem{Tucsnak-Weiss}
M. Tucsnak and G. Weiss, {\em Observation and Control for Operator Semigroups},  Birkh\"auser Verlag, Basel, 2009. 

\bibitem{veraar1}
 J. M. A. M. van Neerven, M. C. Veraar and L. Weis, {\em Stochastic integration in UMD Banach spaces}, Ann. Probab. {\bf 35} (2007), 1438--1478. 
 
\bibitem{Zabczyk}
J. Zabczyk, {\em Bellman's inclusions and excessive measures}, Probab. Math. Statistics {\bf 21} (2001), 101--122. 




\end{thebibliography}
\end{document}